%% file: laphil1.tex
\documentclass[12pt]{article}
\usepackage{amsfonts}
\usepackage{amsmath}
\usepackage{amssymb}

\def\eps{\varepsilon}
\def\R{\mathbb{R}}
\def\N{\mathbb{N}}

\def\Q{\mathbb{Q}}
\def\Z{\mathbb{Z}}
\def\C{\mathbb{C}}

\def\nl{\newline}

\def\bo{\nl\phantom{a}\hfill $\Box$\nl}

\input{psfig.sty}

\begin{document}
\begin{center}
{\Large {\bf Beurling Zeta Functions, Generalised Primes,\\ and
 Fractal Membranes}
 \vspace{0.2in}}\\
{\bf Titus W. Hilberdink and Michel L. Lapidus}
\footnote{The work of M. L. Lapidus was partially supported by the
U. S. National Science Foundation under grant DMS-0070497.}
\end{center}
\indent
\begin{abstract}  
We study generalised prime systems $\mathcal{P}$ 
$(1<p_1\leq p_2\leq\cdots,$ with $p_j\in\R$ tending to infinity) and the associated
Beurling zeta function $\zeta_{\mathcal{P}}(s) =\prod_{j=1}^{\infty} (1-p_j^{-s})^{-1}$. 
Under appropriate assumptions, we establish various analytic properties of 
$\zeta_{\mathcal{P}}(s)$, including its analytic continuation and we characterise the 
existence of a suitable generalised functional equation. In particular, we examine 
the relationship between
a counterpart of the Prime Number Theorem (with error term) and the properties 
of the analytic continuation of $\zeta_{\mathcal{P}}(s)$. Further we study `well-behaved'
g-prime systems, namely, systems for which both the prime and integer counting 
function are asymptotically well-behaved. Finally, we show that there exists a natural
correspondence between generalised prime systems and suitable orders on $\N^2$.

Some of the above results may be relevant to the second author's theory of
`fractal membranes', whose spectral partition functions are precisely given by
Beurling zeta functions.
\vspace{0.3in}\end{abstract}
\begin{center}
{\large {\bf \S 1. Introduction}}
\end{center}
{\bf 1.1 Generalised primes and Beurling zeta functions}\nl
A  {\em generalised prime system} ${\cal P}$ is a sequence of positive reals 
$p_1,p_2,p_3,\ldots$ satisfying
\[ 1<p_1\leq p_2\leq\cdots\leq p_n\leq \cdots\]
and for which $p_n\to\infty$ as $n\to\infty$. From these 
can be formed the system ${\cal N}$ of {\em generalised integers} or {\em Beurling 
integers}; that is, the numbers of the form
\[ p_1^{a_1}p_2^{a_2}\ldots p_k^{a_k}\]
where $k\in\N$ and $a_1,\ldots,a_k\in\N_0$.\footnote{Here and henceforth, $\N =
\{1,2,3,\ldots\}$ and $\N_0=\N\cup\{0\}$.


} For simplicity, we shall often just refer to 
g-primes and g-integers.
 This system generalises the notion of 
prime numbers and the natural numbers obtained from them. Such systems were first
introduced by Beurling \cite{B} and have been studied by numerous authors since then
(see, in particular, the papers by Bateman and Diamond \cite{BD}, Diamond \cite{D1}, 
\cite{D2}, \cite{D3}, \cite{D4}, Hall \cite{H1}, \cite{H2}, Malliavin \cite{M}, Nyman \cite{N}
and Lagarias \cite{Lag}).

Define the {\em counting functions} $\pi_{\cal P}(x)$ and $N_{\cal P}(x)$ by
\[ \pi_{\cal P}(x) = \sum_{p\leq x, p\in\cal P} 1,\tag{1.1}\]
\[ N_{\cal P}(x) = \sum_{n\leq x, n\in\cal N} 1.\tag{1.2}\]
Here, as elsewhere in the paper, we write $\sum_{p\in\cal P}$ to mean a sum over all
the g-primes, counting multiplicities. Similarly for $\sum_{n\in\cal N}$.
Much of the research on this subject has been about connecting the asymptotic 
behaviour of the {\em g-prime counting function} (1.1) and of the 
{\em g-integer counting function} (1.2) as $x\to\infty$. Specifically, given the asymptotic behaviour of 
$\pi_{\cal P}(x)$, what can be said about the behaviour of $N_{\cal P}(x)$, and 
vice versa.

Many of the known results involve the associated {\em zeta function}, 
often referred to as a {\em Beurling zeta function} in the literature, 
which we define formally by the {\em Euler product}
\[ \zeta_{\cal P}(s) = \prod_{p\in\cal P}\frac{1}{1-p^{-s}}.\tag{1.3}\]
This infinite product may be formally multiplied out to give the {\em Dirichlet series}
\[ \zeta_{\cal P}(s) =  \sum_{n\in\cal N} \frac{1}{n^s}.\tag{1.$3^{\prime}$}\]
We are generally interested in those systems for which $\pi_{\cal P}(x) = O(x^A)$ for 
some $A>0$. For then, the series $\sum_p p^{-s}$ converges in a half-plane 
$\{s\in\C:\Re s>\alpha\}$ for some $\alpha\geq 0$. The number $\alpha$ is the abscissa of 
convergence of this series and also of the Euler product for $\zeta_{\cal P}(s)$ in (1.3) 
and of the Dirichlet series $(1.3^{\prime})$. To see 
this, note that for $s$ real and positive, 
\[ \sum_{p\leq x} \frac{1}{p^s}\leq  \sum_{n\leq x} \frac{1}{n^s} \leq 
\exp\biggl\{\sum_{p\leq x}\log\frac{1}{1-p^{-s}}\biggr\}\leq \exp\biggl\{A\sum_{p\leq x}
\frac{1}{p^s}\biggr\}\]
for some positive constant $A$. Hence if one of the series converges as $x\to\infty$, 
then so does the other. We shall further avoid the case $\alpha=0$ for the reason that
it is quite unlike the standard primes and would moreover require quite different methods
to study. Finally, we can renormalise the case $\alpha\in (0,\infty)$ to $\alpha=1$ by
defining a new system ${\cal P}^{\prime} = \{p^{\alpha}:p\in {\cal P}\}$.

We shall usually drop the reference to $\cal P$ and $\cal N$ if there is no chance of
ambiguity. Thus from now on we refer to $\pi(x), N(x), \zeta(s)$ etc. to denote the
functions $\pi_{\cal P}(x), N_{\cal P}(x), \zeta_{\cal P}(s)$.

Note that in the case when $\mathcal{P}$ is the set of (rational) primes, and hence 
$\mathcal{N}$ is the set of natural numbers, then $\zeta(s)$ coincides with the 
classical Riemann zeta function (see, e.g., \cite{E}, \cite{I}, \cite{T2}); further, $\pi(x)$ 
(resp. $N(x)$) is just the standard prime (resp. integer) counting function. 
\nl

Next we briefly discuss the content of the rest of this paper. First we recall 
some of the known results about generalised primes (or integers) and Beurling zeta
functions. We close $\S 1$ by briefly explaining the connections between 
Beurling zeta functions (or g-prime systems) and `fractal membranes'.
\footnote{This last subsection can be omitted on a first reading.}

In $\S 2$, we examine the relationship between a counterpart 
of the Prime Number Theorem in this context and the properties of the analytic
continuation of $\zeta(s)$; we also consider a related question for the integer 
counting function. Further, we study `well-behaved' generalised prime systems,
namely, systems for which (roughly speaking), both the prime and integer counting 
functions are asymptotically well-behaved. 

In $\S 3$, we consider the analytic continuation of $\zeta(s)$, and examine when 
it can be `completed' to satisfy a suitable generalised functional equation. We do not
give a complete answer to the latter difficult question, but indicate several approaches
and give a criterion for the existence of such a functional equation.

Finally, in $\S 4$, we show that there is a natural one-to-one correspondence between 
generalised prime systems and suitable orders on $\N^2$. \nl

\noindent
{\bf 1.2 Summary of relevant known results}\nl
In this section we give a brief synopsis of the known results relating 
the asymptotic behaviours as $x\to\infty$ of $\pi(x)$, $N(x)$, and of some of 
the properties of the Beurling zeta function $\zeta(s)$. 
We start with the connections between $\pi(x)$ and $N(x)$.

The research has concentrated 
on finding conditions for which results of the form
\[ N(x) = ax +E_1(x)\quad\Longleftrightarrow\quad \pi(x) = {\rm li}(x)+E_2(x)\]
hold. Here $a$ is a positive constant, li$(x)$ is the logarithmic integral given by
\[  {\rm li}(x) = \lim_{\eps\to 0^+}\Bigl( \int_0^{1-\eps} + \int_{1+\eps}^x\Bigr)
\frac{dt}{\log t},\]
and $E_1(x)$ and $E_2(x)$ are error terms of smaller order than
$x$ and ${\rm li}(x)$, respectively. The error terms which have been studied (and 
seem to occur naturally) are of three types; namely, those of the form
\[ O\Bigl(\frac{x}{(\log x)^{\gamma}}\Bigr),\,  O(xe^{-c(\log x)^{\alpha}})\quad\mbox{ and }
\quad O(x^{\theta}),\]
where $\gamma>1$, $c>0$ and $\alpha,\theta\in (0,1)$.
\begin{itemize}
\item Beurling (\cite{B}, 1937) showed that
\[ N(x) = ax + O\Bigl(\frac{x}{(\log x)^{\gamma}}\Bigr)\quad\mbox{ for some 
$\gamma>3/2\, \, $ implies}\quad \pi(x) \sim \frac{x}{\log x},
\footnote{Here and henceforth, all such statements are implicitly assumed to be
asymptotic as $x\to\infty$. Moreover, by $f(x)\sim g(x)$, we mean
$f(x)/g(x)\to 1$ as $x\to\infty$.}\]
which is an analogue (in this more general context) of the Prime Number Theorem.
Furthermore, he showed by example that this is false in general for $\gamma = 3/2$.
Conversely, it follows from Diamond's work (\cite{D4}, Theorem 2) that  
\[ \pi(x) = \frac{x}{\log x} + O\Bigl(\frac{x}{(\log x)^{1+\delta}}\Bigr)\quad\mbox{ for some 
$\delta>0\, \, $ implies}\quad N(x) \sim ax.\]
\item Nyman (\cite{N}, 1949) showed that
\[ N(x) = ax + O\Bigl(\frac{x}{(\log x)^A}\Bigr)\quad (\forall A)\quad 
\Longleftrightarrow\quad  \pi(x)  = {\rm li}(x) + O\Bigl(\frac{x}{(\log x)^A}\Bigr)
\quad (\forall A).\]
\item Malliavin (\cite{M}, 1961) showed that
\[ N(x) = ax + O(xe^{-c_1(\log x)^{\alpha}})\tag{1.4}\]
for some $\alpha\in (0,1)$ and $c_1>0$, implies
\[ \pi(x) = {\rm li}(x) + O(xe^{-c_2(\log x)^{\beta}}),\tag{1.5}\]
for some $c_2>0$, where $\beta = \alpha/10$.
Hall (\cite{H2}, 1971) improved this to $\beta = \alpha/7.91$.
Conversely, Malliavin (\cite{M}, 1961) showed that if (1.5) holds for some $\beta\in (0,1)$ 
and $c_2>0$, then (1.4) holds for some $a, c_1>0$ and $\alpha=\frac{\beta}{2+\beta}$.
Diamond (\cite{D2}, 1970) improved this to $\alpha=\frac{\beta}{1+\beta}$, and furthermore,
Diamond's result contains $\log x\log\log x$ instead of $\log x$ in the exponent.  
\item  Landau (\cite{L}, 1903) proved that
\[ N(x) = ax + O(x^{\theta})\qquad\mbox{ for some $\theta<1$}\tag{1.6}\]
implies
\[ \pi(x) = {\rm li}(x) + O(xe^{-c\sqrt{\log x}})\qquad\mbox{ for some $c>0$}.\]
Diamond, Montgomery and Vorhauer have recently shown (see \cite{DMV}) that this is 
essentially best possible by exhibiting a (discrete) system for which (1.6) holds but
\[ \pi(x)-{\rm li}(x) = \Omega(xe^{-c^{\prime}\sqrt{\log x}})\qquad\mbox{ for some 
$c^{\prime}>0$}.\footnote{Here, $f(x)=\Omega(g(x))$ as $x\to\infty$
means that there exists $c>0$ such that $|f(x)|\geq cg(x)$ for some arbitrarily large
$x$.}\]
\end{itemize}

\noindent
{\bf 1.3 Fractal membranes and Beurling zeta functions}\nl
It may be useful to close this introduction by briefly explaining the connections
between the general theme of this paper and the (new) notion of a fractal membrane
(or `quantum fractal string') recently introduced semi-heuristically by the second author
in \cite{Lap-ISRZ} and under current rigorous investigation (by Lapidus and Nest) in
\cite{LapNe1}\footnote{We stress, however, that the present paper 
can be read entirely independently of \cite{Lap-ISRZ} and does not require any
background on the theory of fractal strings or of fractal membranes. It is primarily
a contribution to the theory of Beurling primes and zeta functions.}. Recall that a 
{\em fractal string} $\mathcal{L}=\{l_j\}_{j=1}^{\infty}$ (see e.g. \cite{Lap1-TAMS},
\cite{Lap-Dundee}, \cite{LapPo1}, \cite{LapPo2}, \cite{LapMa}, \cite{HeLap}, 
\cite{Lap-vF2} and \cite{Lap-vF1}, Chapter 1 or \cite{Lap-ISRZ}, \S3.1 for more details 
and motivations) is an open subset of $\R$, 
$\Omega = \cup_{j=1}^{\infty} I_j\subset\R$, whose connected components (i.e.,
the open intervals $I_j$) have lengths $l_j$ such that $\sum_{j=1}^{\infty} l_j^d
<\infty$ for some $d>0$. Without loss of generality,
we may assume that $l_j\downarrow 0$ and $1>l_1\geq l_2\geq\cdots$,
where $\mathcal{L} = \{l_j\}_{j=1}^{\infty}$ is written according to multiplicity.

In the forthcoming book \cite{Lap-ISRZ} --- entitled {\em In Search of the Riemann Zeros}
and building, in particular, upon \cite{Lap-vF1} and \cite{Lap-Springer} --- the 
second author proposes a (physically motivated and noncommutative) geometric
framework within which to `quantize' fractal strings. Given a fractal string $\mathcal{L}$, 
the resulting noncommutative geometric object $\mathcal{T} = \mathcal{T_L}$ --- coined 
a {\em fractal membrane} or {\em quantized fractal string} --- can be thought of 
heuristically as an (adelic, noncommutative) infinite  dimensional torus, with underlying 
`circles' of radii $R_j= 1/\log l_j^{-1}$ ($j=1,2,\ldots$). Further, the spectral (or
quantum) {\em partition function}\footnote{defined as the trace of the `heat semigroup'
associated with $\mathcal{T}$.} of $\mathcal{T}=\mathcal{T}_{\mathcal{L}}$ is given
by
\[ Z_{\mathcal{T}}(s):= \sum_{j=1}^{\infty} e^{-\lambda_j s} = 
\prod_{j=1}^{\infty} (1-l_j^s)^{-1}.\tag{1.7}\]
(See especially, \cite{Lap-ISRZ}, Chapters 3 and 4.) In other words, the spectrum
$\{\lambda_j\}_{j=1}^{\infty}$ of the fractal membrane $\mathcal{T}=
\mathcal{T}_{\mathcal{L}}$ is discrete and is given by the logarithms of the g-integers
$\mathcal{N}=\{n_j\}$ based on the g-prime system $\mathcal{P}=\{p_j\}$, 
with $p_j=l_j^{-1}$ ($j=1,2,\ldots$)\footnote{Physically,
the eigenvalues $\lambda_j$ represent the `energy levels' or the (square root of) the
(vibrational) frequencies of the fractal membrane $\mathcal{T}$.}. With this notation,
$R_j=1/\log p_j$ and the spectral partition function of $\mathcal{T}$ coincides with the
Beurling zeta function of $\mathcal{P}$: $Z_{\mathcal{T}}(s)= \zeta_{\mathcal{P}}(s)$.
Indeed, by (1.7), $Z_{\mathcal{T}}(s)$ is naturally given by the Dirichlet series in 
(1.$3^{\prime}$) and the Euler product in (1.3).

Conjecturally (see \cite{Lap-ISRZ}, \S4.4), $Z_{\mathcal{T}}$ should satisfy a 
generalised functional equation, connecting the (completed) partition functions of
$\mathcal{T}$ and of its `dual fractal membrane' $\mathcal{T}^*$, at the points $s$ and 
$1-s$. In \cite{Lap-ISRZ}, Chapter 5, is also introduced a `moduli space of fractal 
membranes'\footnote{Formally, this is the noncommutative space obtained as the 
quotient of the set of all fractal membranes by the following equivalence relation:
$\mathcal{T}\sim\mathcal{T}'$ iff $\mathcal{P}\sim\mathcal{P}'$ (i.e., $\exists j_0,j_1$ 
such that $p_{j_0+q}=p_{j_1+q}'$, $\forall q\geq 1$). Note that this preserves the poles 
and the zeros (with real part $\neq 0$) of the associated (equivalence class of) partition 
functions.}
which enables one to obtain a natural (noncommutative) flow of zeta functions
(and of g-prime systems), along with the corresponding flow of zeros. Conjecturally,
\cite{Lap-ISRZ}, \S\S 5.4-5.5, this continuous-time `dynamical deformation' of Beurling 
zeta functions and
prime systems would provide a new way to understand the remarkable role played
(within the broader class of Beurling-type zeta functions) by arithmetic (or 
number-theoretic) zeta functions, such as the Riemann zeta function, which 
necessarily satisfy a `self-dual functional equation' (i.e., symbolically, 
$\mathcal{T} =\mathcal{T}^*$).

We mention in closing that in joint work (in preparation) of the second author and
Ryszard Nest, a rigorous operator-algebraic and noncommutative geometric \cite{Con} 
construction of fractal membranes is provided (in \cite{LapNe1})\footnote{It follows from 
\cite{LapNe1} that fractal membranes are truly (second) quantized fractal strings, as 
suggested in \cite{Lap-ISRZ}.} and that (in \cite{LapNe2}) partial results are obtained 
towards some of the conjectures alluded to in the previous paragraph.
From the mathematical point of view, however, much of the program proposed in 
\cite{Lap-ISRZ} remains to be thoroughly investigated.

Although, as was stressed above, they primarily contribute to the standard theory of
Beurling zeta functions and prime systems, the new results obtained in the present 
paper may also contribute to the concrete development of the theory of fractal 
membranes.\nl 

\pagebreak 

\begin{center}
{\large {\bf \S2. Further connections between $\pi(x)$, $N(x)$ and $\zeta(s)$}}
\end{center}

\noindent
{\bf 2.1 Connecting $\pi(x)$ and $\zeta(s)$}\nl
Throughout this paper, we shall use the weighted counting function
\[ \psi_{\cal P}(x) = \sum_{p^k\leq x, k\in\N} \log p = \sum_{n\leq x, n\in\mathcal{N}}
\Lambda(n).\footnote{Here $\Lambda=\Lambda_{\mathcal P}$ denotes the (generalised) 
von Mangoldt function, defined for $n$ in the multiset $\mathcal{N}$ by $\Lambda(n) =
\log p$ if $n=p^m$ for some $p\in\mathcal{P}$ and $m\in\N$, and $\Lambda(n) =0$ 
otherwise.}\]
This is the natural counterpart for a g-prime system $\mathcal{P}$ of the
Chebyshev--von Mangoldt weighted prime power counting function.
As before, we shall often drop the reference to ${\cal P}$ if no confusion is 
likely. In the following, we shall also write
\[ \phi(s) =-\frac{\zeta^{\prime}(s)}{\zeta(s)}= \sum_{n\in\cal N}\frac{\Lambda(n)}{n^s}.\]
The counting functions $N(x)$ and $\psi(x)$ are related to $\zeta(s)$ and $\phi(s)$ via
\[ \zeta(s)=s\int_1^{\infty} \frac{N(x)}{x^{s+1}}\, dx\quad\mbox{ and }\quad
\phi(s)=s\int_1^{\infty} \frac{\psi(x)}{x^{s+1}}\, dx.\]
As a result, it is often more convenient to work with $\psi(x)$, rather than $\pi(x)$. 
Note that for $\alpha\in [\frac{1}{2},1)$, the statements 
\[ \pi(x) = {\rm li}(x) + O(x^{\alpha+\eps})\quad (\forall\,\eps>0)\quad\mbox{ and }\quad
\psi(x) = x + O(x^{\alpha+\eps})\quad (\forall\,\eps>0),\]
are equivalent. For $\mathcal{N}=\N$, it is well-known that the above are equivalent
to the absence of zeros of the Riemann zeta function in the region $\{ s\in \C:\Re s>
\alpha\}$.

We first show that this holds more generally for g-prime systems.\nl

\noindent
{\bf Theorem 2.1}\nl
{\em Suppose that for some $\alpha\in [0,1)$, we have
\[ \psi(x) = x + O(x^{\alpha+\eps})\qquad\mbox{ for all $\eps>0$}.\]
Then $\zeta(s)$ has an analytic continuation to the half-plane
$\{s\in\C:\Re s>\alpha\}$ except for a simple (non-removable) pole at $s=1$ and 
$\zeta(s)\neq 0$ in this region. 

Conversely, suppose that for some $\alpha\in [0,1)$, $\zeta(s)$ has an analytic 
continuation to the half-plane $\{s\in\C:\Re s>\alpha\}$,  except for a simple
(non-removable) pole at $s=1$, and that $\zeta(s)\neq 0$ in this region. Further 
assume that $|\phi(\sigma +it)| = O(|t|^{\eps})$ holds for all $\eps>0$,
uniformly for $\sigma\geq\alpha +\delta$ with any $\delta>0$. Then}
\[ \psi(x) = x + O(x^{\alpha+\eps}) \quad\mbox{{\em for all $\eps>0$}}.\]
The proof is standard except that in the converse part, an extra subtlety arises
due to the possible close proximity of g-integers.\nl

\noindent
{\em Proof.}\,  By hypothesis, $\psi(x) = x+r(x)$ where $r(x) = O(x^{\alpha+\eps})$ for 
all $\eps>0$. It follows that 
\[\phi(s) =s\int_1^{\infty}\frac{x+r(x)}{x^{s+1}}\, dx = \frac{s}{s-1} + 
s\int_1^{\infty}\frac{r(x)}{x^{s+1}}\, dx.\]
The latter integral converges for $\Re s>\alpha$ and represents an analytic function 
in this half-plane. This provides the analytic continuation of
$\phi(s)$ to $\{s\in\C:\Re s>\alpha\}$ except for a simple pole at $s=1$ with residue 1. 
By standard complex analysis, it follows that
%
$\zeta(s)$ has an analytic continuation to $\{ s\in\C:\Re s>\alpha\}\setminus\{1\}$,
except for a simple (non-removable) pole at 1.
Moreover, it has no zeros in this region, for if it did, then
$\phi(s) = \frac{\zeta^{\prime}(s)}{\zeta(s)}$ would have a singularity.


For the converse, note that the hypothesis implies that 
$\phi(s) = -\frac{\zeta^{\prime}(s)}{\zeta(s)}$ has an analytic continuation to 
$\{s\in\C:\Re s>\alpha\}$ except for a simple pole at $s=1$ with residue 1.

Let $c>1$, $T,x>0$ such that $x\not\in\mathcal{N}$. Then, for $n\in\mathcal{N}$,
\[ \frac{1}{2\pi i}\int_{c-iT}^{c+iT} \Bigl(\frac{x}{n}\Bigr)^s\, \frac{ds}{s} = 
O\Bigl(\frac{(x/n)^c}{T|\log x/n|}\Bigr) + \left\{ \begin{array}{cl} 1 & \mbox{ if $n<x$}\\
0&\mbox{ if $n>x$}\end{array}\right. ,\]
where the implied constant is independent of $n$ and $x$.
Multiply through by $\Lambda(n)$ and sum over all $n\in\mathcal{N}$. Thus for $x\not\in
\mathcal{N}$, we have
\[ \frac{1}{2\pi i}\int_{c-iT}^{c+iT} \frac{\phi(s)x^s}{s}\, ds = \psi(x) + O\biggl(\frac{x^c}{T}
\sum_{n\in\mathcal{N}}\frac{\Lambda(n)}{n^c|\log x/n|}\biggr).\]
For $n\leq \frac{x}{2}$ and $n\geq 2x$, $|\log x/n|\geq\log 2$, so
\begin{align*}
 \psi(x) &= \frac{1}{2\pi i}\int_{c-iT}^{c+iT} \frac{\phi(s)x^s}{s}\, ds+O\biggl(\frac{x^c}{T}
\sum_{n\in\mathcal{N}}\frac{\Lambda(n)}{n^c}\biggr)+ O\biggl(\frac{x^c}{T}
\sum_{\frac{x}{2}<n<2x}\frac{\Lambda(n)}{n^c|\log x/n|}\biggr)\\
&= \frac{1}{2\pi i}\int_{c-iT}^{c+iT} \frac{\phi(s)x^s}{s}\, ds+O\Bigl(\frac{x^c}{T(c-1)}\Bigr)+ 
O\biggl(\frac{x\log x}{T}\sum_{\frac{x}{2}<n<2x}\frac{1}{|n- x|}\biggr),\tag{2.1}
\end{align*}
since $\phi(c)=O(\frac{1}{c-1})$ and $|\log x/n|=|\log (1+\frac{n-x}{x})|\asymp
\frac{|n-x|}{x}$ for $\frac{x}{2}<n<2x$. 

Now consider the integral on the right of (2.1). We can push the contour past the pole 
at $s=1$ to the line $\Re s=\sigma$ for any $\sigma>\alpha$. The residue at 1 is $x$. 
Hence 
\[\frac{1}{2\pi i}\int_{c-iT}^{c+iT} \frac{\phi(s)x^s}{s}\, ds = x+ 
\frac{1}{2\pi i}\biggl( \int_{c-iT}^{\sigma-iT}+\int_{\sigma-iT}^{\sigma+iT}+
\int_{\sigma+iT}^{c+iT}\biggr) \frac{\phi(s)x^s}{s}\, ds.\]
We estimate these integrals in turn, using $\phi(s)=O(|t|^{\eps})$. We have
\[ \left|\frac{1}{2\pi i}\int_{\sigma+iT}^{c+iT} \frac{\phi(s)x^s}{s}\, ds\right|\leq\frac{x^c}{2\pi T}
\int_{\sigma}^c |\phi(y+iT)|\, dy = O(x^c T^{-1+\eps}),\]
and similarly for $\int_{c-iT}^{\sigma -iT}$, while
\[ \left|\frac{1}{2\pi i}\int_{\sigma-iT}^{\sigma+iT} \frac{\phi(s)x^s}{s}\, ds\right|\leq
\frac{x^{\sigma}}{2\pi}\int_{-T}^T \frac{|\phi(\sigma+it)|}{\sqrt{\sigma^2+t^2}}\, dt = 
O(x^{\sigma}T^{\eps}).\]
Now choose $c=1+\frac{1}{\log x}$. Then (2.1) gives
\[ \psi(x) = x+O(x T^{-1+\eps})+O(x^{\sigma}T^{\eps})+O\Bigl(\frac{x\log x}{T}\Bigr)
+O\biggl(\frac{x\log x}{T}\sum_{\frac{x}{2}<n<2x}\frac{1}{|n- x|}\biggr)\tag{2.2}\]
for $x\not\in\mathcal{N}$ and every $\eps>0$. We need to bound the sum
on the right but this is difficult in general as $x$ can be arbitrarily close to a 
g-integer. So let's suppose that $x$ is such that there are no g-integers $n$ with
$|n-x|<\frac{1}{x^2}$; i.e. $(x-\frac{1}{x^2},x+\frac{1}{x^2})\cap\mathcal{N}=
\emptyset$. Then
\[  \sum_{\frac{x}{2}<n<2x}\frac{1}{|n- x|}\leq x^2\sum_{\frac{x}{2}<n<2x} 1 \leq
x^2\, N(2x)= O(x^3).\]
Taking $T=x^4$, (2.2) gives $\psi(x) = x+O(x^{\sigma+\eps})$ for all $\eps>0$. This
holds for all $\sigma>\alpha$ so 
\[\psi(x) = x+O(x^{\alpha+\eps})\]
whenever $x\to\infty$ in such a way that $( x-\frac{1}{x^2},x+\frac{1}{x^2})\cap\mathcal{N}
=\emptyset$.

Now we show that this is sufficient to prove the theorem. More precisely, we show the
following: {\em for all $x$ sufficiently large for which $( x-\frac{1}{x^2},x+\frac{1}{x^2})
\cap\mathcal{N}\neq\emptyset$, $\exists$ $x_1\in (x-3,x)$ and $x_2\in (x,x+3)$ such 
that} 
\[ \Bigl( x_1-\frac{1}{x_1^2},x_1+\frac{1}{x_1^2}\Bigr)\cap\mathcal{N}=\emptyset\quad
\mbox{ {\em and} }\quad\Bigl( x_2-\frac{1}{x_2^2},x_2+\frac{1}{x_2^2}\Bigr)\cap
\mathcal{N}=\emptyset.\tag{2.3}\]
Then the result will follow since $x=x_r+O(1)$ and $\psi(x_r)=x_r+O(x_r^{\alpha+\eps})$ 
(for $r=1,2$), so that
\[ \psi(x)\leq\psi(x_2) = x_2+O(x_2^{\alpha+\eps}) = x+O(x^{\alpha+\eps})\]
and
\[ \psi(x)\geq\psi(x_1) = x_1+O(x_1^{\alpha+\eps}) = x+O(x^{\alpha+\eps}).\]
It remains to prove (2.3).

Suppose, for a contradiction, that there is no such $x_2$. Let $y_n=\sqrt[3]{x^3+9n}$, 
for $n\in\N$. Thus each interval $(y_n-\frac{1}{y_n^2},y_n+\frac{1}{y_n^2})$ contains 
an element of $\mathcal{N}$ whenever $y_n< x+3$; i.e. for $n< x^2+3x+3$.
It is elementary to show that 
\[ y_n+\frac{1}{y_n^2}<y_{n+1}-\frac{1}{y_{n+1}^2},\]
so that these intervals are non-overlapping. This means that $N(x+3)-N(x)\geq 
x^2$. But this is false for all $x$ sufficiently large, as $N(x)=O(x)$.

The existence of $x_1$ is shown in a similar way using the sequence $z_n=
\sqrt[3]{x^3-9n}$, leading to $N(x)-N(x-3)\geq x^2$.\bo

Note that if we had assumed the weaker bound $|\phi(\sigma +it)| = O(|t|)$ for $\sigma>
\alpha$, then one would only obtain
\[ \psi(x) = x + O(x^{\frac{\alpha+1}{2}+\eps}) \quad\mbox{ for all $\eps>0$}.\]

\noindent
{\bf 2.2 Connections between $N(x)$ and $\psi(x)$}\nl
In the following, we consider the effect that the assumption on $\psi(x)$ of Theorem 2.1 
has on the asymptotic behaviour of $N(x)$. This extends Diamond's result \cite{D2} 
(relating (1.5) to (1.4)) to the case $\beta=1$. It would be of interest to know if 
(apart from an improved value of $c$) this is essentially best possible.\nl

\noindent
{\bf Theorem 2.2}\nl
{\em Suppose that $\psi(x) = x+O(x^{\alpha})$ for some $\alpha\in (0,1)$. Then
there exist positive constants $\rho$ and $c$ such that}
\[ N(x) = \rho x + O(xe^{-c\sqrt{\log x\log\log x}}).\]

\noindent
{\em Proof.}\, Let $\psi(x)=x+r(x)$, so that $r(x)=O(x^{\alpha})$. We have already seen 
in the proof of Theorem 2.1 that this 
assumption implies the analytic continuation of $\phi(s)$ to $\Re s>\alpha$ except for
a simple pole at $s=1$ with residue 1, and
that for $\Re s>\alpha$, 
\[ s\int_1^{\infty} \frac{r(x)}{x^{s+1}}\, dx = \phi(s)-\frac{s}{s-1}.\]
Now consider the sum $\sum_{n\leq x}\frac{\Lambda(n)}{n^s}$ for $\Re s>\alpha$,
where $n$ ranges over elements of ${\cal N}$. We have
\begin{align*}
\sum_{n\leq x}\frac{\Lambda(n)}{n^s} &=\int_1^x \frac{d\psi(y)}{y^s} =\frac{\psi(x)}{x^s}+
s\int_1^x \frac{\psi(y)}{y^{s+1}}\, dy \\ 
&=x^{1-s}+\frac{r(x)}{x^s}+s\int_1^x \frac{1}{y^s}\, dx
+s\int_1^x \frac{r(y)}{y^{s+1}}\, dy \tag{2.4}\\
&=\frac{x^{1-s}}{1-s}+\phi(s)+\frac{r(x)}{x^s}-
s\int_x^{\infty} \frac{r(y)}{y^{s+1}}\, dy.
\end{align*}
Thus
\[ \phi(s) = \sum_{n\leq x}\frac{\Lambda(n)}{n^s}-\frac{x^{1-s}}{1-s}-\frac{r(x)}{x^s}+
s\int_x^{\infty} \frac{r(y)}{y^{s+1}}\, dy.\]
Writing $s=\sigma +it$, and using $r(x)=O(x^{\alpha})$, we obtain
\[|\phi(\sigma+it)|\leq \sum_{n\leq x}\frac{\Lambda(n)}{n^{\sigma}}+O\Bigl(
\frac{x^{1-\sigma}}{|t|}\Bigr) +O(|t|x^{\alpha-\sigma}).\]
To estimate the first term on the right, put $t=0$ in (2.4) to give
\begin{align*}
\sum_{n\leq x}\frac{\Lambda(n)}{n^{\sigma}}&=\frac{x+r(x)}{x^{\sigma}}+\sigma
\int_1^x y^{-\sigma}\, dy +\sigma\int_1^x \frac{r(y)}{y^{\sigma+1}}\, dy\\
&= x^{1-\sigma}+\frac{\sigma}{1-\sigma}(x^{1-\sigma}-1)+O(x^{\alpha-\sigma})+O(1)\\
&=\frac{x^{1-\sigma}-1}{1-\sigma}+O(1).
\end{align*}
This also holds for $\sigma=1$, if we interpret the first term on the right as $\log x$.
Moreover, with this interpretation, the above estimate is uniform for 
$\sigma\in [\alpha+\delta, c]$ for any $c>1$ and $\delta>0$. Combining these gives
\[ |\phi(\sigma+it)|\leq \frac{x^{1-\sigma}-1}{1-\sigma}+O(1) + O\Bigl(
\frac{x^{1-\sigma}}{|t|}\Bigr) +O(|t|x^{\alpha-\sigma}).\]
The optimal choice for $x$ occurs when $x^{1-\sigma}$ and $|t|x^{\alpha-\sigma}$
are of the same order. So putting $x=|t|^{\frac{1}{1-\alpha}}$, we obtain
\[|\phi(\sigma+it)|\leq \frac{|t|^{\frac{1-\sigma}{1-\alpha}}-1}{1-\sigma}+O(1) + 
O(|t|^{\frac{1-\sigma}{1-\alpha}}).\tag{2.5}\]
Note that for $\sigma=1$ this is
\[|\phi(1+it)|\leq \frac{1}{1-\alpha}\log |t| +O(1).\]
Now we use these inequalities to obtain bounds for $|\zeta(s)|$. For $\sigma\in
(\alpha,1)$,
\begin{align*} 
\log \zeta(\sigma+it) &= \int_{[\sigma+it, 2+it]}\phi(z)\, dz + \log \zeta(2+it)\\
&=\int_{\sigma}^2 \phi(y+it)\, dy +O(1).
\end{align*}
Taking real parts, we obtain
\[\log |\zeta(\sigma+it)| \leq \int_{\sigma}^2 |\phi(y+it)|\, dy +O(1).\]
Let $|t|\geq 3$ and put $\eps(t)=(1-\alpha)\frac{\log\log |t|}{\log|t|}$. Then 
$\alpha<1-\eps(t)<1$. Letting $\sigma=1-\eps(t)$, we deduce from (2.5) that
\[ \log |\zeta(1-\eps(t)+it)| \leq \int_{1-\eps(t)}^2 \frac{|t|^{\frac{1-y}{1-\alpha}}-1}{1-y}\, dy
+O(1)+O\biggl(\int_{1-\eps(t)}^2 |t|^{\frac{1-y}{1-\alpha}}\, dy\biggr).\]
The latter integral equals
\[ \frac{1-\alpha}{\log |t|} (|t|^{\frac{\eps(t)}{1-\alpha}}-|t|^{-\frac{1}{1-\alpha}})<\frac{1-\alpha}
{\log |t|} |t|^{\frac{\eps(t)}{1-\alpha}} =1-\alpha.\]
Hence 
\begin{align*}
\log|\zeta(1-\eps(t)+it)| &\leq \int_{1-\eps(t)}^2 \frac{|t|^{\frac{1-y}{1-\alpha}}-1}{1-y}\, dy
+O(1)\\
&=\int_0^{\eps(t)} \frac{|t|^{\frac{u}{1-\alpha}}-1}{u}\, du + \int_0^1 
\frac{1-|t|^{-\frac{v}{1-\alpha}}}{v}\, dv +O(1)\\
&=\int_0^{\frac{\eps(t)\log|t|}{1-\alpha}} \frac{e^y-1}{y}\, dy + \int_0^{\frac{\log|t|}{1-\alpha}} 
\frac{1-e^{-x}}{x}\, dx +O(1)\\
&= \int_1^{\log\log |t|}\frac{e^y}{y}\, dy +O(\log\log|t|)\\
&\sim\frac{\log|t|}{\log\log|t|}.
\end{align*} 
Hence,
\[ |\zeta(1-\eps(t)+it)| \leq \exp\Bigl(\frac{2\log|t|}{\log\log|t|}\Bigr),\tag{2.6}\]
for all $|t|$ sufficiently large.
We use this bound to find an approximate formula for $N_1(x) =\int_0^x N(y)dy$, via
the formula
\[N_1(x) = \frac{1}{2\pi i}\int_{c-i\infty}^{c+i\infty}\frac{\zeta(s)}{s(s+1)}x^{s+1}\, ds,\]
which holds for any $c>1$. Pushing the contour past the simple pole at 1 gives
\[N_1(x) = \frac{\rho}{2}x^2 +  \frac{1}{2\pi i}\int_{\gamma}\frac{\zeta(s)}{s(s+1)}x^{s+1}\, ds,\]
where $\gamma$ is the contour $s=1-\eps(t)+it$ for $|t|\geq 3$ and $s=1-\eps(3)+it$ for 
$|t|\leq 3$. This is allowed since
\[ \left|\int_{[c+iT, 1-\eps(T) +iT]} \frac{\zeta(s)}{s(s+1)}x^{s+1}\, ds\right|
= O\Bigl(\frac{x^{c+1}}{T}\Bigr) \to 0 \quad\mbox{ as $T\to\infty$.}\]
Now using the bound (2.6) for $|\zeta(s)|$ on $\gamma$, we obtain
\begin{align*}
 \left|N_1(x)-\frac{\rho}{2}x^2\right| &= O\biggl(\int_3^{\infty}\frac{|\zeta(1-\eps(t)+it)|}
{t^2}x^{2-\eps(t)}\, dt\biggr)+O(x^{2-\eps(3)})\\
&=O\biggl(x^2\int_3^{\infty}\exp\left\{ \frac{2\log t}{\log\log t}-2\log t -(1-\alpha)
\frac{\log\log t \log x}{\log t}\right\}\, dt\biggr)+O(x^{2-\eps(3)})\\
&= O\biggl(x^2\int_3^{\infty}\exp\left\{ -\frac{1}{2}u-(1-\alpha)\frac{\log u \log x}{u}\right\}
\, du\biggr)+O(x^{2-\eps(3)}).\tag{2.7}
\end{align*}
Let $\lambda>0$ and split up the integral into two parts with ranges 
$[3,\lambda\sqrt{\log x\log\log x}]$ and $[\lambda\sqrt{\log x\log\log x},\infty)$.
For $u\leq\lambda\sqrt{\log x\log\log x}$, we have
\[ \frac{\log u \log x}{u} \geq \frac{1}{2\lambda}\sqrt{\log x\log\log x},\]
since $\log u/u$ decreases with $u$ for $u\geq e$. Hence
\begin{align*}
\int_3^{\infty}&\exp\left\{ -\frac{1}{2}u-(1-\alpha)\frac{\log u \log x}{u}\right\}
\, du \\
&\leq e^{-\frac{1-\alpha}{2\lambda}\sqrt{\log x\log \log x}} 
\int_3^{\lambda\sqrt{\log x\log\log x}}e^{-\frac{u}{2}}\, du 
+\int_{\lambda\sqrt{\log x\log\log x}}^{\infty}e^{-\frac{u}{2}}\, du\\
&= O( e^{-\frac{1-\alpha}{2\lambda}\sqrt{\log x\log \log x}})+
O( e^{-\frac{\lambda}{2}\sqrt{\log x\log \log x}})\\
&=O( e^{-\lambda^{\prime}\sqrt{\log x\log \log x}}),
\end{align*}
for some $\lambda^{\prime}>0$. In fact, the optimal choice is obtained by taking $\lambda$ 
such that $\lambda=\frac{1-\alpha}{\lambda}$; i.e. $\lambda=\sqrt{1-\alpha}$, which gives 
$\lambda^{\prime}=\frac{1}{2}\sqrt{1-\alpha}$. Hence (2.7) becomes
\[ N_1(x) = \frac{\rho}{2} x^2 +O(x^2 e^{-\frac{1}{2}\sqrt{(1-\alpha)\log x\log \log x}}).\]
By standard methods (using the fact that $N(x)$ increases with $x$), this yields
\[ N(x) = \rho x +O(x e^{-\frac{1}{4}\sqrt{(1-\alpha)\log x\log \log x}}).\]
\bo 

\noindent
{\bf 2.3 Well-behaved systems}\nl
From the comments at the end of \S1.2, we have seen that $N(x)=\rho x+
O(x^{\theta})$ for some $\theta<1$ does not imply that $\psi(x)=x+ O(x^{\theta'})$ for 
some $\theta'<1$. Nor does the converse seem to hold as Theorem 2.2 indicates.
This suggests that we investigate g-prime systems where the functions $N(x)$ and 
$\psi(x)$ are {\em simultaneously} `well-behaved'; that is, for some $\alpha<1$,
\[\psi(x)=x+O(x^{\alpha})\quad\mbox{ and }\quad N(x)=\rho x+O(x^{\alpha}).\]
More precisely, for $0\leq\alpha,\beta<1$, we define an $[\alpha,\beta]$-{\em system} to 
be a g-prime system for which
\begin{align*}
\psi(x) &=x+O(x^{\alpha+\eps})\tag{2.8}\\
N(x) &=\rho x+O(x^{\beta+\eps})\quad \mbox{(for some $\rho>0$)}\tag{2.9}
\end{align*}
hold for all $\eps>0$, but for no $\eps<0$.  
\nl

\noindent
{\bf Conditional Examples}
\begin{enumerate} 
\item For $\N$, (2.9) holds with $\beta=0$ and if the Riemann Hypothesis
were true, (2.8) would hold for $\alpha=\frac{1}{2}$ and this would show the existence of
a $[\frac{1}{2},0]$-system. \nl
\item For the Gaussian integers of the field $\Q(i)$, the Dedekind zeta function is given by
\[ \frac{1}{1-2^{-s}}\prod_{p}\Bigl(\frac{1}{1-p^{-s}}\Bigr)^2\prod_{q}\Bigl(\frac{1}{1-q^{-2s}}
\Bigr)=\frac{1}{4}\sum_{n=1}^{\infty}\frac{r(n)}{n^s},\]
where $p$ and $q$ run over the rational primes 1(mod 4) and 3(mod 4) respectively, and
$r(n)$ is the number of ways of writing $n$ as $a^2+b^2$ with $a,b\in\Z$.
The corresponding prime system $\mathcal{P}$ therefore consists of 2, the rational primes 
$p\equiv 1($mod 4) occurring with multiplicity two, and the squares of the primes 
of the form 3(mod 4). Thus
\[ \pi_{\mathcal P}(x) = 1+2\pi_{1,4}(x)+\pi_{3,4}(\sqrt{x}),\]
where $\pi_{k,m}(x)$ is the number of primes less than or equal to $x$ of the form 
$k$(mod $m)$. On the Generalised Riemann Hypothesis, one has
\[ \pi_{\mathcal P}(x) ={\rm li}(x)+O(x^{\frac{1}{2}+\eps})\mbox{ for all $\eps>0$}.\]
On the other hand, it is known that (see \cite{Hux}) 
\[N_{\mathcal P}(x) = \frac{1}{4}\sum_{n\leq x}r(n) = \frac{\pi}{4}x+O(x^{\frac{23}{73}}),\]
and it is conjectured that the exponent in the error is actually $\frac{1}{4}+\eps$ for all
$\eps>0$. Hence, assuming these conjectures, ${\mathcal P}$ is an example of a 
$[\frac{1}{2},\frac{1}{4}]$-system.
\end{enumerate}

From Theorem 2.1 we see that (2.8) is equivalent to $\zeta(s)$ having a 
simple pole at $s=1$ and having {\em no zeros} in some vertical strip to the left of the 
line $\Re s=1$. Condition (2.9) is equivalent to $\zeta(s)$ being of {\em finite order} in 
such a strip; that is, there exists a positive constant $A$ such that
\[ \zeta(\sigma +it) = O(|t|^A)\quad\mbox{ as $|t|\to\infty$ for $\sigma>\beta$ }.\]
However, if {\em both} (2.8) and (2.9) hold, then the bounds on 
$\zeta(\sigma+it)$ are much stronger, in much the same way as the Riemann 
Hypothesis implies the Lindel\"{o}f Hypothesis.\nl

\noindent
{\bf Theorem 2.3}\nl
{\em Let $\mathcal{P}$ be a $[\alpha,\beta]$-system. Then for 
$\sigma>\Theta=\max\{\alpha,\beta\}$, and uniformly for $\sigma\geq\Theta+\delta$
(any $\delta>0$),}
\[ \phi(\sigma +it) = O((\log |t|)^{\frac{1-\sigma}{1-\Theta}+\eps})\quad\mbox{ and }
\quad \zeta(\sigma +it) = O(\exp\{(\log |t|)^{\frac{1-\sigma}{1-\Theta}+\eps}\}),\]
{\em for all $\eps>0$. In particular, $\zeta(\sigma +it) = O(|t|^{\eps})$ for all $\eps>0$}.\nl

\noindent
{\em Proof.}\, First, $\zeta(s)$ and $\phi(s)$ have analytic continuations to $\{s\in\C:
\Re s>\min\{\alpha,\beta\}\}\setminus\{1\}$ with simple poles at $s=1$. 
Let $s=\sigma+it$, with $\sigma >\min\{\alpha,\beta\}$. Then 
\begin{align*}
\zeta(\sigma+it) &=O(|t|)\quad\mbox{ for $\sigma>\beta$, and}\\
\phi(\sigma+it) &=O(|t|)\quad\mbox{ for $\sigma>\alpha$.}
\end{align*}
Note that for any $\delta>0$, these hold uniformly as $|t|\to\infty$ for 
$\sigma\geq\beta+\delta$, and $\sigma\geq\alpha+\delta$ respectively.

By Theorem 2.1, $\zeta(s)$ is non-zero for $\Re s>\alpha$, and so $\log\zeta(s)$ 
exists and is analytic on $\{s\in\C:\Re s>\alpha\}\setminus (\alpha,1]$.
Hence for $\sigma>\Theta=\max\{\alpha,\beta\}$,
and uniformly for $\sigma\geq\Theta+\delta$, we have
\[ \Re \log \zeta(\sigma+it) = \log |\zeta(\sigma+it)| \leq A\log |t|,\]
for some $A$.
Applying the Borel--Carath\'{e}odory Theorem (see \cite{T1}), it follows that
$|\log \zeta(\sigma+it)| =O(\log |t|)$ uniformly for $\sigma\geq\Theta+2\delta$.
Now by Cauchy's Theorem, 
\[ \phi(\sigma+it) =-\frac{1}{2\pi i}\int_{\gamma}\frac{\log\zeta(z)}{(z-\sigma-it)^2}\, dz,\]
where $\gamma$ is the circle with centre $\sigma+it$ and radius $\eps$. 
Choosing $\eps>0$ so that $\sigma-\eps>\Theta$, gives
\[ |\phi(\sigma+it)|\leq \frac{1}{\eps}\sup_{z\in\gamma}|\log \zeta( z)| =O(\log |t|).\tag{2.10}\]
Let $\delta,\eta>0$.
We apply Hadamard's Three-Circles Theorem to the circles $C_1,C_2,C_3$ with centre
$c+it$ ($c>1+\eta$), passing through the points $1+\eta+it$, $\sigma+it$, $\Theta+\delta +it$
respectively. The radii are $c-1-\eta, c-\sigma,c-\Theta-\delta$ respectively.
Let $M_1,M_2,M_3$ be the maxima of $|\phi(s)|$ on each of the circles $C_1,C_2,C_3$.
Then, 
\[ M_2\leq M_1^{1-\kappa}M_3^{\kappa}, \qquad\mbox{ where }\qquad \kappa=
\frac{\log(\frac{c-\sigma}{c-1-\eta})}{\log(\frac{c-\Theta-\delta}{c-1-\eta})}.\]
Now $M_3 =O(\log |t|)$ by estimate (2.10), and
\[ M_1=\max_{z\in C_1}|\phi(z)|\leq\max_{z\in C_1}\phi(\Re z)\leq \phi(1+\eta)=O(1).\]
Hence $M_2=O((\log|t|)^{\kappa})$, and in particular,
\[ |\phi(\sigma+it)| =O((\log|t|)^{\kappa}).\]
The exponent, $\kappa$, can be made as close as we please to 
$\frac{1-\sigma}{1-\Theta}$ by choosing $c$ large and $\eta,\delta$ small, since
\begin{align*}
 \kappa &=\frac{\log(\frac{c-\sigma}{c-1-\eta})}{\log(\frac{c-\Theta-\delta}{c-1-\eta})}   = 
\frac{\log(\frac{1-\sigma/c}{1-(1+\eta)/c})}{\log(\frac{1-(\Theta+\delta)/c}{1-(1+\eta)/c})}
 = \frac{1+\eta-\sigma+O(\frac{1}{c})}{1+\eta-(\Theta+\delta)+O(\frac{1}{c})}\\
&=\frac{1-\sigma}{1-\Theta}+O(\eta)+O(\delta)+O(1/c).
\end{align*}
Hence
\[ |\phi(\sigma+it)| =O((\log |t|)^{\frac{1-\sigma}{1-\Theta}+\eps})\quad\mbox{ for any 
$\eps>0$}.\]
Finally, 
\begin{align*}
\log|\zeta(\sigma+it)| &=\Re\left\{\int_{[\sigma+it,2+it]}\phi(z)\, dz\right\}+ \log|\zeta(2+it)|
\leq \int_{\sigma}^2 |\phi(x+it)|\, dx+A\\
&=O\Bigl(\int_{\sigma}^2(\log |t|)^{\frac{1-x}
{1-\Theta}+\eps}\, dx\Bigr) =O((\log |t|)^{\frac{1-\sigma}{1-\Theta}+\eps}),
\end{align*}
so that $|\zeta(\sigma+it)| = O(\exp((\log |t|)^{\frac{1-\sigma}{1-\Theta}+\eps}))$. Since
the exponent is less than one for $\sigma>\Theta$ (by taking $\eps$ sufficiently
small), it follows that $|\zeta(\sigma+it)| = O(|t|^{\eps})$ for all $\eps>0$.\bo

\noindent
{\em Remark}\, (i) If $\alpha<\beta$ and we already know that 
$\zeta(s)$ is of finite order for $\sigma>\eta$ for some $\eta\in
(\alpha,\beta)$, then $\zeta(s)$ and $\phi(s)$ have zero 
order in this range.

(ii) If $\beta<\alpha$ and we already know that $\phi(s)$ has 
only finitely many poles for $\sigma>\eta^{\prime}$ (equivalently, 
$\zeta(s)$ has finitely many zeros here), then $\zeta(s)$ 
and $\phi(s)$ have zero order in this range.\nl

%
%
%

For functions $f$ of finite order we define, as usual, the {\em order} $\mu_f(\sigma)$
to be the infimum of all real numbers $\lambda$ such that
\[ f(\sigma+it) = O(|t|^{\lambda})\quad\mbox{ as $|t|\to\infty$.}\]
If there is no such $\lambda$, we shall write $\mu_f(\sigma)=\infty$ and say that $f$ is
of infinite order.
It is well-known that, as a function of $\sigma$, $\mu_f(\sigma)$ is {\em non-negative,
decreasing, and convex}.
Theorem 2.3 tells us that $\mu_{\zeta}(\sigma)=0$ for $\sigma >\Theta$.

We are naturally led to consider the cases where $\alpha=\beta$, $\alpha>\beta$, and 
$\alpha<\beta$. The case $\alpha=\beta$ just tells us that $\zeta(s)$ has no zeros and is
of zero order for $\Re s>\alpha$. 
Of course, we can say nothing about what happens for
$\Re s\leq\alpha$. This leaves us with the more interesting possibilities:\nl

\noindent
{\bf Case A:}\, $\alpha >\beta$.\nl
This looks very similar to the case of $\N$ under the assumption of the Riemann 
Hypothesis. In the half-plane where $\Re s=\sigma>\alpha$, $\mu_{\zeta}(\sigma)=0$
and $\zeta(s)$ has no zeros, while in the strip $\beta<\sigma\leq\alpha$, 
we have $\mu_{\zeta}(\sigma)<1$ and there must be infinitely many zeros on, or 
arbitrarily close to, the line $\Re s=\alpha$. (If not, then, after Remark (ii), in a strip to 
the left of $\Re s=\alpha$, $\zeta(s)$ would be non-zero and of zero order, leading to 
$\psi(x)=x+O(x^{\lambda})$ for some $\lambda<\alpha$ using Theorem 2.1, 
contradicting the minimality of $\alpha$.) See Case A of Figure 1 for a summary of this
discussion.\nl

\noindent
{\bf Case B:}\, $\alpha <\beta$. \nl
This case is quite different. Now $\zeta(s)$ has no zeros for $\Re s>\alpha$, has
zero order for $\Re s>\beta$, and must be of {\em infinite} order in the strip 
$\alpha<\Re s<\beta$. For if $\zeta(s)$ is of finite order in some strip to the left of 
$\beta$, then, by Remark (i), $\mu_{\zeta}(\sigma)=0$ in such a strip. Hence either the 
order is zero or infinite for $\alpha<\Re s<\beta$. However, we can rule out the case of 
zero order since it would imply from (a slightly adjusted) Theorem 2.1 that $N(x) = 
x+O(x^{\lambda})$ for some $\lambda<\beta$, contradicting the minimality of $\beta$.
See Case B of Figure 1.\nl

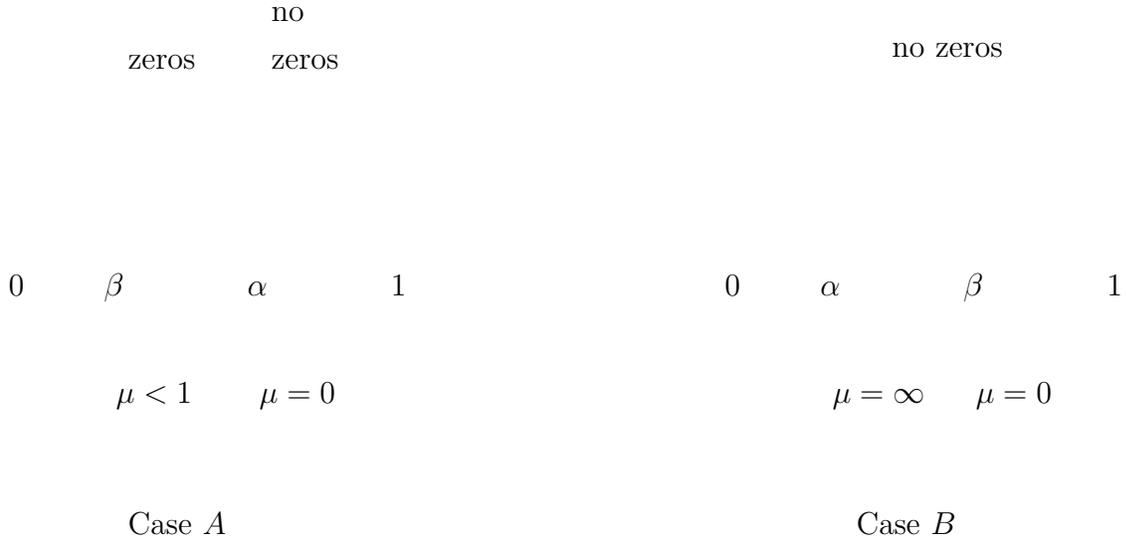
\begin{figure}[h]
\begin{center}
\input{figure2.pstex_t}
\end{center}
\caption{Cases A ($\alpha>\beta$) and B ($\alpha<\beta$) for well-behaved systems}
\label{}
\end{figure}

All the naturally occurring examples of generalised prime systems, such as those 
arising from the Dedekind zeta function, or more generally, from the Selberg class
(\cite{S} and e.g. \cite{CG}, \cite{RM})
of zeta functions (with a standard Euler product), possess a zeta function of finite order. 
Thus assuming the Generalised Riemann Hypothesis, or Selberg's conjecture regarding
the zeros of the corresponding zeta function, we find that all these are examples
of systems $[\frac{1}{2},\beta]$ with $\beta\leq\frac{1}{2}$. In fact, it may reasonably be 
conjectured that these are actually systems $[\frac{1}{2},\frac{1}{2}-\frac{1}{2d}]$, 
where $d$ is the degree of the field extension, or more generally, the 
degree of the Selberg zeta function. 
\nl

\begin{center}
{\large {\bf \S 3. Analytic continuation and Functional equation}}
\end{center}
{\bf 3.1 Analytic continuation}\nl
The classical Riemann zeta function can be continued analytically to the whole plane 
except for a simple pole at $s=1$, so it is natural to ask if, or under what conditions,
this occurs for $\zeta(s)=\zeta_{\mathcal P}(s)$, the Beurling zeta function associated 
to a prime system ${\mathcal P}$. In general, one would expect the 
line $\Re s=1$ to be a natural boundary.

We have already seen (in \S 2) that an analytic continuation exists if either of 
\[ N(x) = \rho x +O(x^{\alpha}), \quad \psi(x) = x +O(x^{\alpha})\]
holds for some $\alpha<1$. 
One problem with such an approach is that there is no hope of extending $\zeta(s)$ 
beyond $\Re s=0$ since $\alpha\geq 0$ in any case.

There are other methods for obtaining analytic continuations. One way is to consider
the function defined for $\Re z>0$ by the series (akin to a `partition function')
\[ F(z) = \sum_{n\in {\cal N}} e^{-nz}.\]
Then
\[ \Gamma(s)\zeta(s) = \sum_{n\in {\cal N}} \int_0^{\infty} x^{s-1}e^{-nx}\, dx =
\int_0^{\infty} x^{s-1}F(x) \, dx\qquad\mbox{ for $\Re s>1$}.\]
The integral $\int_1^{\infty} x^{s-1}F(x) dx$ extends to an entire function since 
$F(x)=O(e^{-x})$, while knowing the behaviour of $F(x)$ near $x=0$ can lead to an 
extension for $\int_0^1 x^{s-1}F(x) dx$. 

Indeed, this method may be generalised by using kernels other than $e^{-x}$; that is, let
\[ F(x) = \sum_{n\in {\cal N}} g(nx),\tag{3.1}\]
for some function $g(x)$ defined on the positive reals with suitable behaviour at infinity. 
By defining the corresponding `gamma' function via the Mellin transform
\[ G(s) = \int_0^{\infty} x^{s-1}g(x)\, dx,\tag{3.2}\]
one finds (formally) that
\[G(s)\zeta(s) = \int_0^{\infty} x^{s-1}F(x)\, dx.\]
We shall suppose $F(x)$ has an asymptotic expansion as $x\to 0^+$ of the form
\[ F(x) \approx \sum_{n=1}^{\infty} \frac{P_n(\log x)}{x^{\lambda_n}},\tag{3.3}\]
where the $P_n(\cdot)$ are polynomials and $\lambda_n$ is a sequence of distinct
complex numbers with real parts decreasing to minus infinity. 
Then we have the following proposition of which we include the proof for completeness.
\nl

\noindent
{\bf Proposition 3.1}\nl
{\em Let $g:(0,\infty)\to\C$ be a continuous function such that $g(x)=O(x^{\alpha})$ as
$x\to 0^+$ and $g(x) = O(x^{-\beta})$ as $x\to\infty$ for some $\alpha<1<\beta$. Let
$F(x)$ and $G(s)$ be defined by $(3.1)$ and $(3.2)$ for $x>0$ and $\Re s\in (\alpha,
\beta)$ respectively. Then
\[ G(s)\zeta(s) = \int_0^{\infty} x^{s-1}F(x)\, dx\qquad\mbox{ for $1<\Re s<\beta$.}\]
Furthermore, if $F(x)$ possesses an asymptotic expansion of the form
$(3.3)$ as $x\to 0^+$, then $G(s)\zeta(s)$ has a continuation to $\{s\in\C:\Re s<\beta\}$ 
which is analytic except for poles at each of the $\lambda_n$, the order of the pole 
being equal to the degree of $P_n(\cdot)$ plus 1.}\nl

\noindent
{\em Proof.}\, 
The integral (3.2) for $G(s)$ converges absolutely and uniformly on compact subsets of
the strip $\alpha<\Re s<\beta$ and $G(s)$ is analytic there, while the series (3.1) for 
$F(x)$ converges absolutely for $x>0$. Thus, we have for $\Re s\in (1,\beta)$,
\[ G(s)\zeta(s) = \sum_{n\in {\cal N}}\frac{1}{n^s}\int_0^{\infty} x^{s-1}g(x)\, dx=
\sum_{n\in {\cal N}}\int_0^{\infty} x^{s-1}g(nx)\, dx = \int_0^{\infty} x^{s-1}F(x)\, dx.\]
The last step of interchanging the integral and sum is justified because of the
absolute convergence of the series. To obtain the analytic continuation, we write, 
\[G(s)\zeta(s) =\int_0^1 x^{s-1}F(x)\, dx+\int_1^{\infty} x^{s-1}F(x)\, dx.\]
The latter integral converges for $\Re s<\beta$ since $F(x)\leq \sum_{n\in {\cal N}}
\frac{A}{(nx)^{\beta}} = O(x^{-\beta})$, so we need only consider the former integral.
By (3.3), we may write
\[ F(x) = F_N(x)+E_N(x),\]
for each $N$ and some $k$, where 
\[F_N(x)=\sum_{n=1}^N \frac{P_n(\log x)}{x^{\lambda_n}}\qquad \mbox{ and }\qquad
E_N(x)=O(x^{-\lambda_N}(\log x)^k).\]
Thus
\begin{align*}
\int_0^1 x^{s-1}F(x)\, dx & = \int_0^1 x^{s-1}F_N(x)\, dx +\int_0^1 x^{s-1}E_N(x)\, dx\\
&=\sum_{n=1}^N \int_0^1 P_n(\log x)x^{s-1-\lambda_n}\, dx +h_N(s),
\end{align*}
where $h_N(s)$ is analytic for $\Re s>\Re \lambda_N$.
But each of the integrals $\int_0^1 P_n(\log x)x^{s-1-\lambda_n}\, dx$ is a polynomial
in $\frac{1}{\lambda_n- s}$ since if $P_n(x) = a_0 + a_1x + \ldots + a_dx^d$, then
\[ \int_0^1 P_n(\log x)x^{s-1-\lambda_n}\, dx = \int_0^{\infty} P_n(-t)e^{-t(s-\lambda_n)}\, 
dt = -\sum_{m=0}^d \frac{m!a_m}{(\lambda_n-s)^{m+1}}. \]
Hence $\int_0^1 x^{s-1}F(x)\, dx$ has an analytic continuation to $\Re s>\lambda_N$
except for poles at $\lambda_1,\ldots,\lambda_N$, whose order is the degree of
$P_n(x)$ plus 1. This holds for all integers $N\geq 1$, and since $\lambda_N\to -\infty$
as $N\to\infty$, we deduce that the meromorphic continuation of $G(s)\zeta(s)$ 
for $\Re s<\beta$. \bo

\input{laphil2.tex}

%% file: figure2.pstex_t
\begin{picture}(0,0)%
\special{psfile=figure2.pstex}%
\end{picture}%
\setlength{\unitlength}{0.012500in}%
\begingroup\makeatletter\ifx\SetFigFont\undefined
\def\x#1#2#3#4#5#6#7\relax{\def\x{#1#2#3#4#5#6}}%
\expandafter\x\fmtname xxxxxx\relax \def\y{splain}%
\ifx\x\y   
\gdef\SetFigFont#1#2#3{%
  \ifnum #1<17\tiny\else \ifnum #1<20\small\else
  \ifnum #1<24\normalsize\else \ifnum #1<29\large\else
  \ifnum #1<34\Large\else \ifnum #1<41\LARGE\else
     \huge\fi\fi\fi\fi\fi\fi
  \csname #3\endcsname}%
\else
\gdef\SetFigFont#1#2#3{\begingroup
  \count@#1\relax \ifnum 25<\count@\count@25\fi
  \def\x{\endgroup\@setsize\SetFigFont{#2pt}}%
  \expandafter\x
    \csname \romannumeral\the\count@ pt\expandafter\endcsname
    \csname @\romannumeral\the\count@ pt\endcsname
  \csname #3\endcsname}%
\fi
\fi\endgroup
\begin{picture}(495,315)(75,445)
\put(245,545){\makebox(0,0)[lb]{\smash{\SetFigFont{12}{14.4}{rm}1}}}
\put(545,545){\makebox(0,0)[lb]{\smash{\SetFigFont{12}{14.4}{rm}1}}}
\put(130,500){\makebox(0,0)[lb]{\smash{\SetFigFont{12}{14.4}{rm}$\mu<1$}}}
\put(190,500){\makebox(0,0)[lb]{\smash{\SetFigFont{12}{14.4}{rm}$\mu =0$}}}
\put(440,445){\makebox(0,0)[lb]{\smash{\SetFigFont{12}{14.4}{rm}${\rm Case}\,\, B$}}}
\put(430,500){\makebox(0,0)[lb]{\smash{\SetFigFont{12}{14.4}{rm}$\mu=\infty$}}}
\put(490,500){\makebox(0,0)[lb]{\smash{\SetFigFont{12}{14.4}{rm}$\mu=0$}}}
\put(135,445){\makebox(0,0)[lb]{\smash{\SetFigFont{12}{14.4}{rm}${\rm Case}\,\, A$}}}
\put( 85,545){\makebox(0,0)[lb]{\smash{\SetFigFont{12}{14.4}{rm}$0$}}}
\put(135,640){\makebox(0,0)[lb]{\smash{\SetFigFont{12}{14.4}{rm}zeros}}}
\put(455,645){\makebox(0,0)[lb]{\smash{\SetFigFont{12}{14.4}{rm}no zeros}}}
\put(125,545){\makebox(0,0)[lb]{\smash{\SetFigFont{12}{14.4}{rm}$\beta$}}}
\put(185,545){\makebox(0,0)[lb]{\smash{\SetFigFont{12}{14.4}{rm}$\alpha$}}}
\put(425,545){\makebox(0,0)[lb]{\smash{\SetFigFont{12}{14.4}{rm}$\alpha$}}}
\put(485,545){\makebox(0,0)[lb]{\smash{\SetFigFont{12}{14.4}{rm}$\beta$}}}
\put(385,545){\makebox(0,0)[lb]{\smash{\SetFigFont{12}{14.4}{rm}$0$}}}
\put(195,660){\makebox(0,0)[lb]{\smash{\SetFigFont{12}{14.4}{rm}no}}}
\put(195,640){\makebox(0,0)[lb]{\smash{\SetFigFont{12}{14.4}{rm}zeros}}}
\end{picture}

%% file: laphil2.tex
\noindent
{\bf 3.2 Functional Equation for $\zeta(s)$}\nl
Riemann's $\zeta(s)$ satisfies the well-known functional equation
\[ \zeta(1-s) = 2^{1-s}\pi^{-s}\cos\Bigl(\frac{\pi s}{2}\Bigr)\Gamma(s)\zeta(s)\]
which may also be written in the form
\[ \pi^{-(1-s)/2}\Gamma\Bigl(\frac{1-s}{2}\Bigr)\zeta(1-s)=\pi^{-s/2}
\Gamma\Bigl(\frac{s}{2}\Bigr)\zeta(s).\tag{3.4}\]
We can ask whether some form
of functional equation still holds for more general Beurling zeta functions 
$\zeta_{\mathcal P}(s)$. Specifically, under what
circumstances does an equation of the form
\[G_1(1-s)\zeta_1(1-s)=G_2(s)\zeta_2(s),\tag{3.5}\]
hold? Here the $G_r(s)$ are `Gamma'-like functions, say defined by Mellin transforms
of given functions $g_r(x)$ as in (3.2), and the $\zeta_r(s)$ are Euler products for two 
(possibly different) prime systems. It is well-known that (3.4) is equivalent to the modular 
identity:
\[ \sum_{-\infty}^{\infty} e^{-\pi n^2x^2} = \frac{1}{x}\sum_{-\infty}^{\infty} e^{-\pi n^2/x^2}.\]
It follows from Theorem 3.2 below that, under some mild conditions on the $g_r(x)$, 
(3.5) is also equivalent to an identity of this type. For similar equivalences in a 
related context, see \cite{Ber}.\nl

\noindent
{\bf Theorem 3.2}\nl
{\em For $r=1,2$, let $F_r:(0,\infty)\to\C$ be continuous functions which are 
$O(x^{-1-\eps})$ as $x\to 0^+$ for every $\eps>0$ and $O(x^{-\beta_r})$ as $x\to\infty$, 
for some $\beta_r>1$. Let $\Psi_r(s)$ be the Mellin transform of $F_r(x)$ defined for 
$\Re s\in (1,\beta_r)$. Then the following are equivalent:}

\begin{enumerate}
\item {\em $\Psi_1(s)$ and $\Psi_2(s)$ have analytic continuations to 
$1-\beta_2<\Re s<\beta_1$ and $1-\beta_1<\Re s<\beta_2$ respectively except for a 
finite number of poles in any given strip $\sigma_1\leq\Re s\leq\sigma_2$, tend to
$0$ uniformly as $|\Im s|\to\infty$, and satisfy the functional equation}
\[\Psi_1(1-s)=\Psi_2(s)\qquad\mbox{ {\em for} $1-\beta_1<\Re s<\beta_2$}.\]
\item
\[ F_1(x) = \frac{1}{x}F_2\Bigl(\frac{1}{x}\Bigr) + H(x),\tag{3.6}\]
{\em where $H(x)$ is a finite series of the form $\sum_k a_kx^{\mu_k} (\log x)^{\nu_k}$ 
with $\mu_k\in\C$ and $\nu_k\in\N_0$.}
\end{enumerate}
For the specific form of functional equation (3.5), we take 
\[F_r(x) = \sum_{n\in\mathcal{N}_r} g_r(nx)\qquad (r=1,2),\]
where $\mathcal{N}_r$ are two g-integer systems, the $g_r(x)$ are as before and satisfy 
$O(x^{-\alpha_r})$ and $O(x^{\beta_r})$ as $x\to 0^+$ and $x\to\infty$ respectively, for 
some $\alpha_r<1<\beta_r$. Then $F_r(x)$ has the appropriate behaviour at $0$ and 
infinity, and $\Psi_r(s)=G_r(s)\zeta_r(s)$. \nl

{\em Addendum}\,
After the completion of this paper, it was pointed out to us that S. Bochner (\cite{Bo},
Theorems 2 and 3) had obtained a very similar result. There is a minor difference.
In Bochner's case, the function $\Psi_1(s)$ is allowed to have infinitely many poles in a given strip rather than finitely many. This implies that $H(x)$ does not have the closed form that we obtain, but is a more general `residual' function.\nl

\noindent
{\em Proof of Theorem 3.2.}\,  (a)$\Longrightarrow$(b).\, 
By the inverse Mellin transform for $F_2(x)$, we have
\[ F_2(x) = \frac{1}{2\pi i}\int_{(c)}\Psi_2(s)x^{-s}\, ds\qquad\mbox{ for any $c\in(1,
\beta_2)$,}\]
where $\int_{(c)}$ denotes $\lim_{T\to\infty} \int_{c-iT}^{c+iT}$. Hence
\begin{align*}
\frac{1}{x}F_2\Bigl(\frac{1}{x}\Bigr) &= \frac{1}{2\pi i}\int_{(c)}\Psi_2(s)x^{s-1}\, ds\\
=& \frac{1}{2\pi i}\int_{(c)}\Psi_1(1-s)x^{s-1}\, ds\\
=& \frac{1}{2\pi i}\int_{(c^{\prime})}\Psi_1(s)x^{-s}\, ds\qquad\mbox{(where $c^{\prime}
=1-c)$.}
\end{align*}
On the other hand,
\[ F_1(x) = \frac{1}{2\pi i}\int_{(c)}\Psi_1(s)x^{-s}\, ds\qquad\mbox{ for any $c\in(1,
\beta_1)$.}\]
Moving the contour from the line $\Re s=c$ to $\Re s=c^{\prime}$, we pick up the 
residues at the poles of  $\Psi_1(s)$ in the strip $\{ s:0\leq \Re s\leq 1\}$. 
By assumption there are only a finite number of these and the 
residues of the integrand are all of the form $ax^{\mu}P(\log x)$ for some $a,
\mu\in\C$ and some polynomial $P(\cdot)$. Moving the contour is permissible here 
by the assumption that $\Psi_1(\sigma +iT)\to 0$ uniformly as $|T|\to\infty$. Hence
\[F_1(x) = \frac{1}{x}F_2\Bigl(\frac{1}{x}\Bigr) + \sum_k a_kx^{\mu_k} (\log x)^{\nu_k}\]
for some constants $a_k,\mu_k\in\C$ and $\nu_k\in\N_0$. \nl

\noindent
(b)$\Longrightarrow$(a).\, For $r=1,2$, we write
\[ \Psi_r(s) =  \int_0^1 x^{s-1}F_r(x)\, dx +\int_1^{\infty} x^{s-1}F_r(x)\, dx = A_r(s)+B_r(s),\]
where $A_r(s) = \int_0^1 x^{s-1}F_r(x)\, dx$ and $B_r(s) = \int_1^{\infty} x^{s-1}F_r(x)\, dx$.
These integrals converge for $\Re s>1$ and $\Re s<\beta_r$ respectively.
Now 
\begin{align*}
A_1(s) &= \int_0^1 x^{s-1}F_1(x)\, dx = \int_0^1 x^{s-1}\Bigl(\frac{1}{x}F_2\Bigl(\frac{1}{x}
\Bigr)+H(x)\Bigr)\, dx\\
&= \int_1^{\infty} x^{-s}F_2(x)\, dx +\int_0^1 x^{s-1}H(x)\, dx\\
& = B_2(1-s)+H_1(s).\tag{3.7}
\end{align*}
On the other hand,
\begin{align*}
B_1(s) &= \int_1^{\infty} x^{s-1}F_1(x)\, dx = \int_1^{\infty} x^{s-1}\Bigl(\frac{1}{x}
F_2\Bigl(\frac{1}{x}\Bigr)+H(x)\Bigr)\, dx\\
&= \int_0^1 x^{-s}F_2(x)\, dx +\int_1^{\infty} x^{s-1}H(x)\, dx\\
& = A_2(1-s)+H_2(s).\tag{3.8}
\end{align*}
Where do (3.7) and (3.8) hold? The functions $H_1(s)$ and $H_2(s)$ are actually
rational functions as shown below and hence have meromorphic continuations to $\C$.
Now $A_1(s)$ is analytic for $\Re s>1$ while $B_2(1-s)$ is analytic for $\Re s>1-
\beta_2$. Thus $A_1(s)$ has a meromorphic continuation for $\Re s>1-\beta_2$ and
(3.7) holds for $\Re s>1-\beta_2$. Similarly, (3.8) holds for $\Re s<\beta_1$. 
These provide the required continuations of $\Psi_1(s)$ and $\Psi_2(s)$. In particular,
for $1-\beta_2 <\Re s<\beta_1$, (3.7) and (3.8) hold simultaneously and together give
\[ A_1(s) + B_1(s) =  A_2(1-s)+B_2(1-s)+H_1(s)+H_2(s),\]
i.e.
\[ \Psi_1(s) = \Psi_2(1-s) +H_1(s)+H_2(s).\]
We proceed to show that $H_1(s)$ and $H_2(s)$ are rational functions such that 
$H_1(s)+H_2(s)=0$. Indeed, we have
\begin{align*} H_1(s) &= \int_0^1 x^{s-1}H(x)\, dx = \sum_k a_k\int_0^1 x^{s+\mu_k -1}
(\log x)^{\nu_k}\, dx \tag{$x=e^{-t}$}\\
&=\sum_k a_k(-1)^{\nu_k}\int_0^{\infty} t^{\nu_k}e^{-(s+\mu_k)t}\, dt,
\end{align*}
the integral converging whenever $\Re (s+\mu_k)>0$ for all $k$. For $H_2(s)$ we 
have 
\begin{align*} H_2(s) &= \int_1^{\infty} x^{s-1}H(x)\, dx = \sum_k a_k\int_1^{\infty}
x^{s+\mu_k -1}(\log x)^{\nu_k}\, dx \tag{$x=e^t$}\\
&=\sum_k a_k\int_0^{\infty} t^{\nu_k}e^{(s+\mu_k)t}\, dt,
\end{align*}
the integral converging whenever $\Re (s+\mu_k)<0$ for all $k$. But for $m$ a 
non-negative integer and $\alpha\in\C$ with positive real part, $\int_0^{\infty} t^m 
e^{-\alpha t}dt =m!\alpha^{-m-1}$. Thus
\[H_1(s) = \sum_k \frac{a_k(-1)^{\nu_k}\nu_k!}{(s+\mu_k)^{\nu_k+1}},\qquad
\mbox{ and}\]
\[H_2(s) = \sum_k \frac{a_k\nu_k!}{(-(s+\mu_k))^{\nu_k+1}} = \sum_k 
\frac{a_k(-1)^{\nu_k+1}\nu_k!}{(s+\mu_k)^{\nu_k+1}}= -H_1(s),\]
as required.

Finally, we show that $|\Psi_r(s)|\to 0$ uniformly as $|\Im s|\to\infty$, for $r=1,2$. 
By the functional 
equation, it suffices to show this for $|\Psi_1(s)|$. Let $s=\sigma +it$. We have
\[ \Psi_1(s)= A_1(s)+B_1(s)= B_2(1-s)+B_1(s)+H_1(s).\]
Now $|H_1(s)| = O(1/|t|)$, while 
\[ |B_1(\sigma+it)| = \left|\int_1^{\infty} x^{\sigma-1+it}F_1(x)\, dx\right| =
\left|\int_1^{\infty} x^{\sigma-1}F_1(x)e^{it\log x}\, dx\right|.\]
The integral $\int_1^{\infty} x^{\sigma-1}F_1(x)dx$ converges uniformly whenever 
$\sigma\leq \beta_1-\delta$ for any $\delta>0$; hence by the Riemann--Lebesgue 
Theorem, it follows that the RHS integral tends to 0 as $|t|\to\infty$ uniformly for 
$\sigma\leq \beta_1-\delta$. Similarly, $|B_2(1-\sigma+it)|\to 0$ uniformly for 
$\sigma\geq 1-\beta_2+\delta$. Hence $|\Psi_1(\sigma+it)|\to 0$ uniformly as $|t|\to
\infty$ for $1-\beta_2+\delta\leq\sigma\leq \beta_1-\delta$.
\bo

\noindent
{\em Remark.}\, In a work in preparation \cite{LapNe2}, the second author and R. Nest  
have obtained (under hypotheses similar to those of Proposition 3.1) a generalised 
functional equation for $\zeta(s)$. The latter connects a suitable completion of the 
Beurling zeta function $\zeta(s)$ associated to the g-prime system $\mathcal{P}$ 
to that of $\zeta^*(s)$, associated to a `dual system' $\mathcal{P}^*$ (which involves 
in general a {\em continuous} g-prime system);
\footnote{Determining when these g-prime systems have an {\em 
infinite discrete part} is an interesting and difficult problem that is far from resolved
at this stage.} this establishes part of a conjecture formulated in \cite{Lap-ISRZ}, 
\S4.4. (see \S 1.3 above).
We note that the context of \cite{LapNe2} is not restricted to Beurling zeta functions
and applies, in particular, to zeta functions associated with quasicrystals. 
\nl

\noindent
{\bf 3.2 Functional Equation for $\zeta(s)$}\nl
Riemann's $\zeta(s)$ satisfies the well-known functional equation
\[ \zeta(1-s) = 2^{1-s}\pi^{-s}\cos\Bigl(\frac{\pi s}{2}\Bigr)\Gamma(s)\zeta(s)\]
which may also be written in the form
\[ \pi^{-(1-s)/2}\Gamma\Bigl(\frac{1-s}{2}\Bigr)\zeta(1-s)=\pi^{-s/2}
\Gamma\Bigl(\frac{s}{2}\Bigr)\zeta(s).\tag{3.4}\]
We can ask whether some form
of functional equation still holds for more general Beurling zeta functions 
$\zeta_{\mathcal P}(s)$. Specifically, under what
circumstances does an equation of the form
\[G_1(1-s)\zeta_1(1-s)=G_2(s)\zeta_2(s),\tag{3.5}\]
hold? Here the $G_r(s)$ are `Gamma'-like functions, say defined by Mellin transforms
of given functions $g_r(x)$ as in (3.2), and the $\zeta_r(s)$ are Euler products for two 
(possibly different) prime systems. It is well-known that (3.4) is equivalent to the modular 
identity:
\[ \sum_{-\infty}^{\infty} e^{-\pi n^2x^2} = \frac{1}{x}\sum_{-\infty}^{\infty} e^{-\pi n^2/x^2}.\]
It follows from Theorem 3.2 below that, under some mild conditions on the $g_r(x)$, 
(3.5) is also equivalent to an identity of this type. For similar equivalences in a 
related context, see \cite{Ber}.\nl

\noindent
{\bf Theorem 3.2}\nl
{\em For $r=1,2$, let $F_r:(0,\infty)\to\C$ be continuous functions which are 
$O(x^{-1-\eps})$ as $x\to 0^+$ for every $\eps>0$ and $O(x^{-\beta_r})$ as $x\to\infty$, 
for some $\beta_r>1$. Let $\Psi_r(s)$ be the Mellin transform of $F_r(x)$ defined for 
$\Re s\in (1,\beta_r)$. Then the following are equivalent:}

\begin{enumerate}
\item {\em $\Psi_1(s)$ and $\Psi_2(s)$ have analytic continuations to 
$1-\beta_2<\Re s<\beta_1$ and $1-\beta_1<\Re s<\beta_2$ respectively except for a 
finite number of poles in any given strip $\sigma_1\leq\Re s\leq\sigma_2$, tend to
$0$ uniformly as $|\Im s|\to\infty$, and satisfy the functional equation}
\[\Psi_1(1-s)=\Psi_2(s)\qquad\mbox{ {\em for} $1-\beta_1<\Re s<\beta_2$}.\]
\item
\[ F_1(x) = \frac{1}{x}F_2\Bigl(\frac{1}{x}\Bigr) + H(x),\tag{3.6}\]
{\em where $H(x)$ is a finite series of the form $\sum_k a_kx^{\mu_k} (\log x)^{\nu_k}$ 
with $\mu_k\in\C$ and $\nu_k\in\N_0$.}
\end{enumerate}
For the specific form of functional equation (3.5), we take 
\[F_r(x) = \sum_{n\in\mathcal{N}_r} g_r(nx)\qquad (r=1,2),\]
where $\mathcal{N}_r$ are two g-integer systems, the $g_r(x)$ are as before and satisfy 
$O(x^{-\alpha_r})$ and $O(x^{\beta_r})$ as $x\to 0^+$ and $x\to\infty$ respectively, for 
some $\alpha_r<1<\beta_r$. Then $F_r(x)$ has the appropriate behaviour at $0$ and 
infinity, and $\Psi_r(s)=G_r(s)\zeta_r(s)$. \nl

\noindent
{\em Proof of Theorem 3.2.}\,  (a)$\Longrightarrow$(b).\, 
By the inverse Mellin transform for $F_2(x)$, we have
\[ F_2(x) = \frac{1}{2\pi i}\int_{(c)}\Psi_2(s)x^{-s}\, ds\qquad\mbox{ for any $c\in(1,
\beta_2)$,}\]
where $\int_{(c)}$ denotes $\lim_{T\to\infty} \int_{c-iT}^{c+iT}$. Hence
\begin{align*}
\frac{1}{x}F_2\Bigl(\frac{1}{x}\Bigr) &= \frac{1}{2\pi i}\int_{(c)}\Psi_2(s)x^{s-1}\, ds\\
=& \frac{1}{2\pi i}\int_{(c)}\Psi_1(1-s)x^{s-1}\, ds\\
=& \frac{1}{2\pi i}\int_{(c^{\prime})}\Psi_1(s)x^{-s}\, ds\qquad\mbox{(where $c^{\prime}
=1-c)$.}
\end{align*}
On the other hand,
\[ F_1(x) = \frac{1}{2\pi i}\int_{(c)}\Psi_1(s)x^{-s}\, ds\qquad\mbox{ for any $c\in(1,
\beta_1)$.}\]
Moving the contour from the line $\Re s=c$ to $\Re s=c^{\prime}$, we pick up the 
residues at the poles of  $\Psi_1(s)$ in the strip $\{ s:0\leq \Re s\leq 1\}$. 
By assumption there are only a finite number of these and the 
residues of the integrand are all of the form $ax^{\mu}P(\log x)$ for some $a,
\mu\in\C$ and some polynomial $P(\cdot)$. Moving the contour is permissible here 
by the assumption that $\Psi_1(\sigma +iT)\to 0$ uniformly as $|T|\to\infty$. Hence
\[F_1(x) = \frac{1}{x}F_2\Bigl(\frac{1}{x}\Bigr) + \sum_k a_kx^{\mu_k} (\log x)^{\nu_k}\]
for some constants $a_k,\mu_k\in\C$ and $\nu_k\in\N_0$. \nl

\noindent
(b)$\Longrightarrow$(a).\, For $r=1,2$, we write
\[ \Psi_r(s) =  \int_0^1 x^{s-1}F_r(x)\, dx +\int_1^{\infty} x^{s-1}F_r(x)\, dx = A_r(s)+B_r(s),\]
where $A_r(s) = \int_0^1 x^{s-1}F_r(x)\, dx$ and $B_r(s) = \int_1^{\infty} x^{s-1}F_r(x)\, dx$.
These integrals converge for $\Re s>1$ and $\Re s<\beta_r$ respectively.
Now 
\begin{align*}
A_1(s) &= \int_0^1 x^{s-1}F_1(x)\, dx = \int_0^1 x^{s-1}\Bigl(\frac{1}{x}F_2\Bigl(\frac{1}{x}
\Bigr)+H(x)\Bigr)\, dx\\
&= \int_1^{\infty} x^{-s}F_2(x)\, dx +\int_0^1 x^{s-1}H(x)\, dx\\
& = B_2(1-s)+H_1(s).\tag{3.7}
\end{align*}
On the other hand,
\begin{align*}
B_1(s) &= \int_1^{\infty} x^{s-1}F_1(x)\, dx = \int_1^{\infty} x^{s-1}\Bigl(\frac{1}{x}
F_2\Bigl(\frac{1}{x}\Bigr)+H(x)\Bigr)\, dx\\
&= \int_0^1 x^{-s}F_2(x)\, dx +\int_1^{\infty} x^{s-1}H(x)\, dx\\
& = A_2(1-s)+H_2(s).\tag{3.8}
\end{align*}
Where do (3.7) and (3.8) hold? The functions $H_1(s)$ and $H_2(s)$ are actually
rational functions as shown below and hence have meromorphic continuations to $\C$.
Now $A_1(s)$ is analytic for $\Re s>1$ while $B_2(1-s)$ is analytic for $\Re s>1-
\beta_2$. Thus $A_1(s)$ has a meromorphic continuation for $\Re s>1-\beta_2$ and
(3.7) holds for $\Re s>1-\beta_2$. Similarly, (3.8) holds for $\Re s<\beta_1$. 
These provide the required continuations of $\Psi_1(s)$ and $\Psi_2(s)$. In particular,
for $1-\beta_2 <\Re s<\beta_1$, (3.7) and (3.8) hold simultaneously and together give
\[ A_1(s) + B_1(s) =  A_2(1-s)+B_2(1-s)+H_1(s)+H_2(s),\]
i.e.
\[ \Psi_1(s) = \Psi_2(1-s) +H_1(s)+H_2(s).\]
We proceed to show that $H_1(s)$ and $H_2(s)$ are rational functions such that 
$H_1(s)+H_2(s)=0$. Indeed, we have
\begin{align*} H_1(s) &= \int_0^1 x^{s-1}H(x)\, dx = \sum_k a_k\int_0^1 x^{s+\mu_k -1}
(\log x)^{\nu_k}\, dx \tag{$x=e^{-t}$}\\
&=\sum_k a_k(-1)^{\nu_k}\int_0^{\infty} t^{\nu_k}e^{-(s+\mu_k)t}\, dt,
\end{align*}
the integral converging whenever $\Re (s+\mu_k)>0$ for all $k$. For $H_2(s)$ we 
have 
\begin{align*} H_2(s) &= \int_1^{\infty} x^{s-1}H(x)\, dx = \sum_k a_k\int_1^{\infty}
x^{s+\mu_k -1}(\log x)^{\nu_k}\, dx \tag{$x=e^t$}\\
&=\sum_k a_k\int_0^{\infty} t^{\nu_k}e^{(s+\mu_k)t}\, dt,
\end{align*}
the integral converging whenever $\Re (s+\mu_k)<0$ for all $k$. But for $m$ a 
non-negative integer and $\alpha\in\C$ with positive real part, $\int_0^{\infty} t^m 
e^{-\alpha t}dt =m!\alpha^{-m-1}$. Thus
\[H_1(s) = \sum_k \frac{a_k(-1)^{\nu_k}\nu_k!}{(s+\mu_k)^{\nu_k+1}},\qquad
\mbox{ and}\]
\[H_2(s) = \sum_k \frac{a_k\nu_k!}{(-(s+\mu_k))^{\nu_k+1}} = \sum_k 
\frac{a_k(-1)^{\nu_k+1}\nu_k!}{(s+\mu_k)^{\nu_k+1}}= -H_1(s),\]
as required.

Finally, we show that $|\Psi_r(s)|\to 0$ uniformly as $|\Im s|\to\infty$, for $r=1,2$. 
By the functional 
equation, it suffices to show this for $|\Psi_1(s)|$. Let $s=\sigma +it$. We have
\[ \Psi_1(s)= A_1(s)+B_1(s)= B_2(1-s)+B_1(s)+H_1(s).\]
Now $|H_1(s)| = O(1/|t|)$, while 
\[ |B_1(\sigma+it)| = \left|\int_1^{\infty} x^{\sigma-1+it}F_1(x)\, dx\right| =
\left|\int_1^{\infty} x^{\sigma-1}F_1(x)e^{it\log x}\, dx\right|.\]
The integral $\int_1^{\infty} x^{\sigma-1}F_1(x)dx$ converges uniformly whenever 
$\sigma\leq \beta_1-\delta$ for any $\delta>0$; hence by the Riemann--Lebesgue 
Theorem, it follows that the RHS integral tends to 0 as $|t|\to\infty$ uniformly for 
$\sigma\leq \beta_1-\delta$. Similarly, $|B_2(1-\sigma+it)|\to 0$ uniformly for 
$\sigma\geq 1-\beta_2+\delta$. Hence $|\Psi_1(\sigma+it)|\to 0$ uniformly as $|t|\to
\infty$ for $1-\beta_2+\delta\leq\sigma\leq \beta_1-\delta$.
\bo

\noindent
{\em Remark.}\, In a work in preparation \cite{LapNe2}, the second author and R. Nest  
have obtained (under hypotheses similar to those of Proposition 3.1) a generalised 
functional equation for $\zeta(s)$. The latter connects a suitable completion of the 
Beurling zeta function $\zeta(s)$ associated to the g-prime system $\mathcal{P}$ 
to that of $\zeta^*(s)$, associated to a `dual system' $\mathcal{P}^*$ (which involves 
in general a {\em continuous} g-prime system);
\footnote{Determining when these g-prime systems have an {\em 
infinite discrete part} is an interesting and difficult problem that is far from resolved
at this stage.} this establishes part of a conjecture formulated in \cite{Lap-ISRZ}, 
\S4.4. (see \S 1.3 above).
We note that the context of \cite{LapNe2} is not restricted to Beurling zeta functions
and applies, in particular, to zeta functions associated with quasicrystals. 
\nl

\begin{center}
{\large {\bf \S 4. Partial orders on $\cal N$}}
\end{center}
We have of course a great deal of freedom when choosing a system of generalised primes. 
Each $p_j$ can be chosen arbitrarily as long as it is larger than $p_{j-1}$. 
Every choice of system ${\cal P}$ results in an {\em ordering} of ${\cal N}$. For two 
different systems, we would (in general) expect to have two different orderings on
the corresponding $\cal N$s. For example, in one case we might have $p_1^2<p_2$, while
in the next maybe $p_1^2>p_2$. This raises the following 
questions:
\begin{enumerate} 
\item Are orderings on $\mathcal{N}$ uniquely determined by the choice of $\cal P$? 
\item Given that we know the ordering on $\cal N$, can we reconstruct $\cal P$?
\end{enumerate}
The answer to (a) is negative for an obvious reason: having chosen $\cal P$, then
for any $\lambda>0$, the system ${\cal P}^{\lambda} = \{p^{\lambda}:p\in {\cal P}\}$  
produces the same ordering in its generalised integers as that of $\cal P$ --- they
are just the $\lambda^{\rm th}$ powers of the integers in $\cal N$. We show below that 
this is the only case where this happens; that is, two essentially different prime
systems produce different orderings (see also \cite{F} for similar results). 
In the following, we set ${\mathcal N}^{\lambda} = \{ n^{\lambda}:n\in\mathcal{N}\}$.\nl

\noindent
{\bf Theorem 4.1}\nl
{\em Let ${\cal P}_1$ and ${\cal P}_2$ be two generalised prime systems with
generalised integers ${\cal N}_1$ and ${\cal N}_2$ respectively. If the orderings of 
${\cal N}_1$ and ${\cal N}_2$ coincide, then ${\cal P}_1= {\cal P}_2^{\lambda}$
for some $\lambda>0$ and hence ${\cal N}_1= {\cal N}_2^{\lambda}$}.\nl

\noindent
{\em Proof.}\, Denote the primes in ${\cal P}_1$ by $p_1, p_2, \ldots$, and those in 
${\cal P}_2$ by $q_1, q_2, \ldots$. Let $k,m\in\N$. Then $p_k^m\in [p_1^n,p_1^{n+1})$ 
for some $n\in\N$. Since ${\cal N}_2$ has the same ordering as ${\cal N}_1$ we also 
have $q_k^m\in [q_1^n,q_1^{n+1})$. Taking logs gives
\[ n\leq\frac{m\log p_k}{\log p_1}<n+1\quad\mbox{ and }\quad n\leq
\frac{m\log q_k}{\log q_1}<n+1.\]
i.e. $n= [\frac{m\log p_k}{\log p_1}] = [\frac{m\log p_k}{\log p_1}]$. It follows that
\[ \frac{m\log p_k}{\log p_1}-1<n\leq \frac{m\log q_k}{\log q_1} <n+1\leq 
\frac{m\log p_k}{\log p_1}+1,\]
and hence
\[ \frac{\log p_k}{\log p_1} - \frac{1}{m}<\frac{\log q_k}{\log q_1}<
\frac{\log p_k}{\log p_1}+\frac{1}{m}.\]
This holds for all $m$. Letting $m\to\infty$ gives
\[ \frac{\log p_k}{\log p_1}=\frac{\log q_k}{\log q_1}, \]
i.e. $p_k=q_k^{\lambda}$ where $\lambda =\frac{\log p_1}{\log q_1}$, and hence
${\cal P}_1= {\cal P}_2^{\lambda}$.\bo

The above theorem actually shows that the order is determined uniquely (up to a 
constant) by the way the {\em powers} of the generalised primes are ordered.
It is therefore of importance to consider the set of powers of generalised primes.
We shall denote it by $\cal Q$; i.e.
\[ {\cal Q} = \{ p^n:p\in {\cal P}, n\in\N\}.\]
The set $\cal Q$ is isomorphic to $\N^2=\N\times\N$ via the isomorphism:
$p_m^n\mapsto (m,n)$. Furthermore, the ordering on $\cal Q$ induces an order
on $\N^2$. Such an order necessarily satisfies the following three axioms:
\begin{description}
\item[A1] $(m,n)\leq (m^{\prime},n^{\prime})$ whenever both $m\leq m^{\prime}$ and
$n\leq n^{\prime}$, with strict inequality if $n<n^{\prime}$.
\item[A2] $(m,n)\leq (m^{\prime},n^{\prime})$ implies $(m,kn)\leq (m^{\prime},
kn^{\prime})$ for every $k\in\N$, with strict inequality if $(m,n)< 
(m^{\prime},n^{\prime})$.
\item[A3] Finiteness: (i) for all $n\in\N$, there exists $k\in\N$ such that $(1,n)<(k,1)$ and\nl
\hspace*{.71in}(ii) for all $m\in\N$, there exists $l\in\N$ such that $(m,1)<(1,l)$.
\end{description}
On the other hand, we show below in Theorem 4.3 that for any order on $\N^2$ 
satisfying axioms A1, A2, and A3, there is a generalised prime system $\cal P$ for 
which $\cal Q$ has the same ordering.

First we need the following lemma:
\nl

\noindent
{\bf Lemma 4.2}\nl
{\em Let $f:\N\to\R$ be such that for all $m,n\in\N$
\[ \left|\frac{f(mn)}{mn}-\frac{f(n)}{n}\right|\leq \frac{1}{n}.\tag{4.1}\]
Then $\lim_{n\to\infty} \frac{f(n)}{n}$ exists.}\nl

\noindent
{\em Proof.}\, By putting $n=1$ in (4.1), it follows that $\frac{f(n)}{n}$ is bounded.
Let
\[\alpha = {\lim\inf}_{n\to\infty}\frac{f(n)}{n}.\]
By definition of $\alpha$, there exists a sequence $n_k$ tending to infinity with $k$ 
such that 
\[\frac{f(n_k)}{n_k}\to\alpha\quad\mbox{ as $k\to\infty$}.\]
Hence, for fixed $m\in\N$, we have $\frac{f(mn_k)}{mn_k}\to\alpha$ as $k\to\infty$, since
\[ \Biggl|\frac{f(mn_k)}{mn_k}-\underbrace{\frac{f(n_k)}{n_k}}_{\to\alpha}\Biggr|\leq 
\underbrace{\frac{1}{n_k}}_{\to 0}.\]
Now fix $n\in\N$, put $m=n_k$ in (4.1) and let $k\to\infty$. Then
\[ \Biggl|\underbrace{\frac{f(nn_k)}{nn_k}}_{\to\alpha}-\frac{f(n)}{n}\Biggr|\leq \frac{1}{n}.\]
Thus
\[ \left|\frac{f(n)}{n}-\alpha\right|\leq \frac{1}{n}\]
and the result follows.\bo

\noindent
{\bf Theorem 4.3}\nl
{\em Given an order on $\N^2$ satisfying axioms A1-A3, there exists a generalised 
prime system $\cal P$ which induces this order.}\nl

\noindent 
{\em Proof.}\, Fix $k\in\N$. For $n\in\N$, let $f_k(n)$ be the unique positive integer for 
which
\[ (1,f_k(n))\leq (k,n)<(1,f_k(n)+1).\]
This exists on account of A2 and A3, and is unique by A1.
Hence, replacing $n$ by $mn$ ($m\in\N$), we have 
\[ (1,f_k(mn))\leq (k,mn)<(1,f_k(mn)+1)\tag{4.2}.\]
On the other hand, A2 implies that
\[  (1,mf_k(n))\leq (k,mn)<(1,m(f_k(n)+1))\tag{4.3}.\]
$(4.2)$ and $(4.3)$ give
\begin{align*}
(1,f_k(mn)) &<(1,m(f_k(n)+1)),\mbox{ and}\\
(1,mf_k(n)) &<(1,f_k(mn)+1).
\end{align*}
Thus $f_k(mn)<mf_k(n)+m$ and $mf_k(n)<f_k(mn)+1$. Combining these and dividing 
through by $mn$ gives
\[ \frac{f_k(n)}{n} - \frac{1}{mn} <\frac{f_k(mn)}{mn} <\frac{f_k(n)}{n} +\frac{1}{n},\tag{4.4}\]
and so $f_k(n)$ satisfies condition (4.1) of Lemma 4.2. Hence $\frac{f_k(n)}{n} \to\alpha_k$
for some $\alpha_k$. Also, $\alpha_k$ increases with $k$, since $f_k(n)$ does.

Choose $p_1>1$ arbitrarily, and define g-primes by $p_k = p_1^{\alpha_k}$
(note that $\alpha_1=1$). This gives a system of generalised primes which induces
the same order. 
Indeed, let $m\to\infty$ in (4.4). Then $f_k(n)\leq n\alpha_k\leq f_k(n)+1$.
Now if $p_m^n <p_{m'}^{n'}$ (i.e. $n\alpha_m < n' \alpha_{m'}$), then 
\[ f_m(n)\leq n\alpha_m < n^{\prime}\alpha_{m'}\leq f_{m'}(n^{\prime})+1.\]
Since, $f_m(n)$ and  $f_{m'}(n^{\prime})$ are integers, this implies $f_m(n)\leq
f_{m'}(n^{\prime})$, i.e. $(m,n)\leq (m^{\prime}, n^{\prime})$.
But if $(m,n)= (m^{\prime}, n^{\prime})$, then $(m,kn)= 
(m^{\prime}, kn^{\prime})$ for all $k\in\N$ by A2. Hence $\frac{f_m(kn)}{k}=
\frac{f_{m'}(kn^{\prime})}{k}$ and, letting $k\to\infty$, we have
$n\alpha_m = n' \alpha_{m'}$, i.e. $p_m^n =p_{m'}^{n'}$. 

This shows that $p_m^n <p_{m'}^{n'}$ implies $(m,n)< (m^{\prime}, n^{\prime})$.
\bo

\noindent
{\em Remark.}\, In the light of the results obtained in the present section, it would 
be interesting to develop a suitable theory of generalised valuations for g-integer
systems. We leave this problem for future investigations.\nl

\noindent
Titus W. Hilberdink, Department of Mathematics, University of Reading, Whiteknights,
 PO Box 220, Reading RG6 6AX, UK.\nl
E-mail address: t.w.hilberdink@reading.ac.uk\nl

\noindent
Michel L. Lapidus, Department of Mathematics, University of California, Riverside,
CA 92521-0135, USA.\nl
E-mail: lapidus@math.ucr.edu

\end{document}

%% file: laphil1.bbl
{\small
\begin{thebibliography}{99}
\bibitem{Apos} T. M. Apostol, {\em Mathematical Analysis}, Addison-Wesley,
Massachusetts, 1973.
\bibitem{BD} P. T. Bateman and H. G. Diamond, Asymptotic distribution of Beurling's
generalised prime numbers, in: {\em Studies in Number Theory} 6, Prentice-Hall, 1969,
pp. 152-212. 
\bibitem{Ber} B. C. Berndt, Identities involving the coefficients of a class of 
Dirichlet series I, {\em Trans. Amer. Math. Soc.} {\bf 137} (1969), 345-359.
\bibitem{B} A. Beurling, Analyse de la loi asymptotique de la distribution des
nombres premiers g\'{e}n\'{e}ralis\'{e}s, I, {\em Acta Math.} {\bf 68} (1937), 255-291.
\bibitem{Bo} S. Bochner, Some properties of modular relations, {\em Ann. of Math.}
{\bf 53} No. 2 (1951), 332-363.
\bibitem{Con} A. Connes, {\em Noncommutative Geometry}, Academic Press, New York,
1994.
\bibitem{CG} J. B. Conrey and A. Ghosh, On the Selberg class of Dirichlet series: small 
degrees, {\em Duke Math. J.} {\bf 72} (1993), 673-693.
\bibitem{D1} H. G. Diamond, The prime number theorem for Beurling generalised 
numbers, {\em J. Number Theory} {\bf 1} (1969), 200-207. 
\bibitem{D2} H. G. Diamond, Asymptotic distribution of Beurling generalised integers,
{\em Illinois J. Math.} {\bf 14} (1970), 12-28.
\bibitem{D3} H. G. Diamond, A set of generalised numbers showing Beurling's
theorem to be sharp, {\em Illinois J. Math.} {\bf 14} (1970), 29-34. 
\bibitem{D4} H. G. Diamond, When do Beurling generalised integers have a 
density?, {\em J. Reine Angew. Math.} {\bf 295} (1977), 22-39.
\bibitem{DMV} H. Diamond, H. Montgomery and U. Vorhauer, Beurling primes with large 
oscillation. Preprint, 2003.
\bibitem{E} H. M. Edwards, {\em Riemann's Zeta Function}, Academic Press, 
New York, 1974. (Reprinted by Dover Publ., 2001.)
\bibitem{F} L. Fuchs, {\em Partially Ordered Algebraic Systems}, International 
Series of Monographs in Pure and Applied Mathematics 28, Pergamon Press, 1963. 
\bibitem{H1} R. S. Hall, The prime number theorem for generalised primes, 
{\em J. Number Theory} {\bf 4} (1972), 313-320.
\bibitem{H2} R. S. Hall, Theorems about Beurling's generalised primes and the 
associated zeta function, Ph.D. Thesis, University of Illinois, Urbana, 1967.
\bibitem{HeLap} C. Q. He and M. L. Lapidus, Generalized
Minkowski content, spectrum of fractal drums, fractal
strings and the Riemann zeta function, {\em Memoirs Amer.
Math. Soc.}, No. 608, {\bf 127} (1997), 1-97.
\bibitem{H} T. W. Hilberdink, Some connections between Bernoulli convolutions and
analytic number theory, in: {\em Fractal Geometry and Number Theory: A Jubilee of 
Benoit Mandelbrot} (M. L. Lapidus and M. van Frankenhuysen, eds.), Proc. Sympos.
Pure Math., Amer. Math. Soc., Providence R. I., 2004 (to appear).
\bibitem{Hux} M. N. Huxley, Exponential sums and lattice points, {\em Proc. London 
Math. Soc.} {\bf 60} (1990), 471-475.
\bibitem{I} A. E. Ingham, {\em The Distribution of Prime Numbers}, Cambridge Tracts in 
Math. 30, Second edition (reprinted from the 1932 edition), Cambridge University Press, 
1990.
\bibitem{K} J. Knopfmacher, {\em Abstract Analytic Number Theory}, Dover, 1990.
\bibitem{Lag} J. C. Lagarias, Beurling Generalised Integers with the Delone Property,
{\em Forum Math.} {\bf 11} (1999), 295-312.
\bibitem{L} E. Landau, Neuer Beweis des Primzahlsatzes und Beweis des
Primidealsatzes, {\em Math. Ann.} {\bf 56} (1903), 645-670.
\bibitem{Lap1-TAMS} M. L. Lapidus, Fractal drums, inverse spectral problems for elliptic 
operators and a partial resolution of the Weyl--Berry conjecture, {\em Trans. Amer. 
Math. Soc.} {\bf 325} (1991), 465-529.
\bibitem{Lap-Dundee} M. L. Lapidus, Vibrations of fractal drums, the Riemann
hypothesis, waves in fractal media, and the Weyl--Berry conjecture, in: {\em Ordinary
and Partial Differential Equations} (B. D. Sleeman and R. J. Jarvis, eds.), vol. IV,
Proc. Twelfth Internat. Conf. (Dundee, Scotland, UK, June 1992), Pitman Research
Notes in Math. Series, vol. 289, Longman Scientific and Technical, London, 1993, 
pp. 126-209.
\bibitem{Lap-Springer} M. L. Lapidus, T-duality, functional equation, and 
noncommutative stringy spacetime, in: {\em New Interactions of Mathematics with the
Natural Sciences and the Humanities} (L. Boi, ed.), Springer-Verlag, Berlin, 2003,
pp. 1-92 (in press).
\bibitem{Lap-ISRZ} M. L. Lapidus, {\em In Search of the Riemann Zeros: Strings, fractal
membranes, and noncommutative spacetimes}, Book in preparation, 2003, 365+(xii) pp.
\bibitem{LapMa} M. L. Lapidus and H. Maier, The Riemann hypothesis and inverse 
spectral problems for fractal strings, {\em J. London Math. Soc.} (2) {\bf 52} (1995), 15-34.
\bibitem{LapNe1} M. L. Lapidus and R. Nest, Fractal membranes and the second 
quantization of fractal strings, in preparation. 
\bibitem{LapNe2} M. L. Lapidus and R. Nest, Functional equations for zeta functions 
associated with quasicrystals and fractal membranes, in preparation.
\bibitem{LapPo1} M. L. Lapidus and C. Pomerance, The Riemann zeta-function and the 
one-dimensional Weyl--Berry conjecture for fractal drums, {\em Proc. London Math. Soc.}
(3) {\bf 66} (1993), 41-69.
\bibitem{LapPo2} M. L. Lapidus and C. Pomerance, Counterexamples to the 
modified Weyl--Berry conjecture on fractal drums, {\em Proc. Cambridge Philos. Soc.}
{\bf 199} (1996), 167-178.
\bibitem{Lap-vF1} M. L. Lapidus and M. van Frankenhuysen, {\em Fractal
Geometry and Number Theory: Complex dimensions of fractal strings and zeros of
zeta functions}, Birkh\"{a}user, Boston, 2000. (Second revised and enlarged edition
to appear in 2004.)
\bibitem{Lap-vF2} M. L. Lapidus and M. van Frankenhuysen, Fractality, self-similarity
and complex dimensions, in: {\em Fractal Geometry and Number Theory: A Jubilee 
of Benoit Mandelbrot} (M. L. Lapidus and M. van Frankenhuysen, eds.), Proc. Sympos.
Pure Math., Amer. Math. Soc., Providence, R. I., 2004 (to appear).
\bibitem{M} P. Malliavin, Sur le reste de la loi asymptotique de r\'{e}partition des
nombres premiers g\'{e}n\'{e}ralis\'{e}s de Beurling, {\em Acta Math.} {\bf 106} (1961),
281-298.
\bibitem{N} B. Nyman, A general prime number theorem, {\em Acta Math.} {\bf 81} (1949),
299-307.
\bibitem{RM} M. Ram Murty, Selberg's conjectures and Artin $L$-functions, 
{\em Bull. Amer. Math. Soc.} {\bf 31} (1994), 1-14.
\bibitem{R} C. Ryavec, The analytic continuation of Euler products with applications
to asymptotic formulae, {\em Illinois J. Math.} {\bf 17} (1973), 608-616.
\bibitem{S} A. Selberg, Old and new conjectures about a class of Dirichlet series, in: 
A. Selberg: {\em Collected Papers}, vol. II, Springer-Verlag, New York, 1991, pp. 47-63.
\bibitem{T1} E. C. Titchmarsh, {\em The Theory of Functions}, Second edition,
Oxford University Press, 1986.
\bibitem{T2} E. C. Titchmarsh, {\em The Theory of the Riemann Zeta-function},
Second edition, Oxford University Press, 1986.
\bibitem{Z1} Wen-Bin Zhang, Density and $O$-Density of Beurling generalised
integers, {\em J. Number Theory} {\bf 30} (1988), 120-139.
\bibitem{Z2} Wen-Bin Zhang, A generalization of Hal\'{a}sz's theorem to Beurling's
generalised integers and its application, {\em Illinois J. Math.} {\bf 31}
(1987), 645-664.
\end{thebibliography}
}
